\documentclass[11pt,english]{article}
\usepackage[T1]{fontenc}
\usepackage[utf8]{inputenc}
\usepackage{geometry}
\usepackage{color}
\definecolor{note_fontcolor}{rgb}{0.800781, 0.800781, 0.800781}
\usepackage{babel}
\usepackage{array}
\usepackage{units}
\usepackage{mathrsfs}
\usepackage{multirow}
\usepackage{amsmath}
\usepackage{amsthm}
\usepackage{accents}
\usepackage{ytableau}
\usepackage{amssymb}
\usepackage{bbm}
\usepackage{relsize}
\usepackage{stackrel}
\usepackage[all]{xy}
\PassOptionsToPackage{normalem}{ulem}
\usepackage{ulem}
\usepackage[unicode=true,pdfusetitle,
 bookmarks=true,bookmarksnumbered=false,bookmarksopen=false,
 breaklinks=false,pdfborder={0 0 1},backref=page,colorlinks=true]
 {hyperref}
\hypersetup{
 linkcolor=blue,citecolor=blue}

\makeatletter

\newtheorem{theorem}{Theorem}[section]
\newtheorem{cor}[theorem]{Corollary}
\newtheorem{lem}[theorem]{Lemma}
\newtheorem{prop}[theorem]{Proposition}
\newtheorem{remark}[theorem]{Remark}
\newtheorem{defn}[theorem]{Definition}

\newtheorem{example}[theorem]{Example}

\newtheorem{claim}[theorem]{Claim}

\numberwithin{equation}{section}
\numberwithin{figure}{section}

\usepackage{babel}
\usepackage{ytableau}

\usepackage{multirow}

\usepackage{appendix}

\makeatother

\begin{document}
\global\long\def\F{\mathbb{\mathbb{\mathbf{F}}}}%

\global\long\def\cfq{{\cal CF}_q}%
\global\long\def\scfq{{\cal SCF}_q}%

\global\long\def\scf{{\cal SCF}_\textrm{sym}}%
 
\global\long\def\rk{\mathbb{\mathrm{rk}}}%
 
\global\long\def\Hom{\mathrm{Hom}}%
 
\global\long\def\defi{\stackrel{\mathrm{def}}{=}}%
 
\global\long\def\tr{{\cal T}r }%
 
\global\long\def\id{\mathrm{id}}%
 
\global\long\def\Aut{\mathrm{Aut}}%
 
\global\long\def\wl{w_{1},\ldots,w_{\ell}}%
 
\global\long\def\alg{\le_{\mathrm{alg}}}%
 
\global\long\def\A{{\cal A}}%
 
%\global\long\def\mobius{M\dacute{o}bius}%
 
\global\long\def\chimax{\chi^{\mathrm{max}}}%
 
\global\long\def\uexp{\mathbb{E}_{\mathrm{unif}}}%
 
\global\long\def\symirr{\widehat{S_{\infty}}}%
 
\global\long\def\supp{\mathrm{supp}}%
  
\global\long\def\suppi{\widetilde{\mathrm{supp}}}%
 
\global\long\def\fq{\mathbb{F}_{q}}%
 
\global\long\def\IT{I|_{T}}%
 
\global\long\def\gl{\mathrm{GL}}%
 
\global\long\def\glm{\mathrm{GL}_{m}\left(q\right)}%
 
\global\long\def\gln{\mathrm{GL}_{N}\left(q\right)}%
 
\global\long\def\fq{\mathbb{F}_q}
 
\global\long\def\fix{\mathrm{fix}}%

\global\long\def\M{\mathrm{Mat}_{m \times N}\left(\fq \right)}%
 
\global\long\def\R{{\cal R}}
\global\long\def\P{{\cal P}}
\global\long\def\I{{\cal I}}
 
\global\long\def\indep{\mathfrak{indep}}%
 
\global\long\def\B{B}%
 
\global\long\def\btil{\tilde{\B}}%

\global\long\def\etil{\tilde{{e}}}%
 
\global\long\def\free{\mathfrak{Free}}%
 
\global\long\def\C{{\cal C}}%
 
\global\long\def\lm{\mathfrak{lm}}%
 
\global\long\def\spn{\mathrm{span}}%

\global\long\def\chark{\mathrm{char}_k\left(\B \right)}%
 
\global\long\def\sw{{\cal S}_{w}}%
 
\global\long\def\fr{\mathrm{f.r.}}%

\global\long\def\rkk{\mathrm{rank(k)}}%

\global\long\def\rkm{\mathrm{rank(m)}}%

\global\long\def\triv{\mathrm{triv}}%

\global\long\def\chimu{\chi_{\vec{\mu}[\bullet]}}%

\title{The ring of stable characters over $\protect\gl_{\bullet}\left(q\right)$}
\author{Danielle Ernst-West~~~~~~Doron Puder~~~~~~Yotam Shomroni}
\maketitle

%%%%%%%%%%%%%%%%%%%%%%%%%%%%%%%%%%%%%%%%%%%%%%%%%%%%%%%%%%%%%
%
%  Abstract
%
%%%%%%%%%%%%%%%%%%%%%%%%%%%%%%%%%%%%%%%%%%%%%%%%%%%%%%%%%%%%%

\begin{abstract}
For a fixed prime power $q$, let $\gl_{\bullet}\left(q\right)$ denote the family of groups $\gl_N\left(q\right)$ for $N \in \mathbb{Z}_{\geq 0}$. In this paper we study the $\mathbb{C}$-algebra of "stable" class functions of $\gl_{\bullet}\left(q \right)$, and show it admits four different linear bases, each arising naturally in different settings. One such basis is that of stable irreducible characters, namely, the class functions spanned by the characters corresponding to finitely generated simple $\mathrm{VI}$-modules in the sense of \cite{putman2017representation,gan2018representation}. A second one comes from characters of parabolic representations. The final two, one originally defined in \cite{gadish2020dimension} and the other in  \cite{ernstwest2022word}, are more combinatorial in nature. As corollaries, we clarify many properties of these four bases and prove a conjecture from \cite{balachandran2023product}.

\end{abstract}

%\tableofcontents

%%%%%%%%%%%%%%%%%%%%%%%%%%%%%%%%%%%%%%%%%%%%%%%%%%%%%%%%%%%%%
%
%  1. Introduction
%
%%%%%%%%%%%%%%%%%%%%%%%%%%%%%%%%%%%%%%%%%%%%%%%%%%%%%%%%%%%%%

\section{Introduction\label{sec:intro}}
Throughout this paper we fix a prime power $q$ and consider the family of general linear groups $\gl_N\left(q\right)$ over a fixed finite field $\mathbb{F}_q$ of order $q$ for $N \in \mathbb{Z}_{\geq 0}$, where $\gl_0(q)$ is the trivial group.
Denote by $\gl_{\bullet}(q)$ the union of all $\gln$ for $N \in \mathbb{Z}_{\geq 0},$ namely, $\gl_{\bullet}(q) = \bigcup_{N \in \mathbb{Z}_{\geq 0}} \gl_{N}(q)$. For brevity, we write $G_N$ for the group $\gln$ for $N\in\mathbb{Z}_{\ge0}$, and $G_\bullet$ for $\gl_{\bullet}(q)$.

This paper studies complex-valued functions which are defined simultaneously on $G_N$ for all but finitely many values of $N$, and that are, moreover, class functions on $G_N$ whenever they are defined. Formally, our functions will belong to the $\mathbb{C}$-algebra $\cfq$ defined as follows:

\begin{defn}[${\cal CF}_q$] \label{def:CF_q}
    Let $f=\{f_N\}_{N\ge N_0}$ where $N_0\in \mathbb{Z}_{\ge0}$  and $f_N\colon G_N\to\mathbb{C}$ for all $N\ge N_0$. If $f_N$ is a class function on $G_N$ for all $N\ge N_0$, we say that $f$ is a class function on $G_\bullet$. We say that two class functions $f$ and $g$ on $G_\bullet$ are equivalent if $f_N=g_N$ for every large enough $N$. 
    
    The equivalence classes of class function form a $\mathbb{C}$-algebra with respect to point-wise addition and multiplication. We denote this algebra by $\cfq$.
\end{defn}

For brevity of notation, we shall often treat class functions $f$ on $G_\bullet$ as if they were elements of $\cfq$, even though the actual element of $\cfq$ is the equivalence class $[f]$. 

We now define four different subalgebras of $\cfq$: $\R$, $\R^\fr$, $\P$ and $\I$, which all turn out to be identical by our main result -- Theorem \ref{thm:main}.

%%%%%%%%%%%%%%%%%%%%%%%%%%%%%%%%%%%%%%%%%%%%%%%%%%%%%%%%%%%%%
% R
%%%%%%%%%%%%%%%%%%%%%%%%%%%%%%%%%%%%%%%%%%%%%%%%%%%%%%%%%%%%%
\subsection*{The algebra $\R$ and its subspaces $\R_m$ and $\R_{\le m}$}

In \cite{ernstwest2022word}, the first two authors together with Seidel studied word measures on $G_N$: these are probability measures on these groups defined by word maps. In particular, various families of real- or complex-valued functions
defined on $G_N$ were considered, and their expected value under word measures was studied.
The main family of functions studied in [Ibid.] are the functions counting the number of $m$-tuples of vectors in  $V_N\stackrel{\text{def}}{=}\mathbb{F}_q^{~N}$ which span a subspace which is invariant under $g$ and on which $g$ acts in
a prescribed way. This is formalized as follows:

\begin{defn}
\label{def:B_w} For every  $m\in \mathbb{Z}_{\ge0}$ and $\B\in G_m$, define a map 
$\btil_N\colon G_N \to\mathbb{Z}_{\ge0}$, for arbitrary $N\in \mathbb{Z}_{\ge0}$ as follows. For $g\in G_N$ we let $\btil_N\left(g\right)$
be the number of $m$-tuples of row vectors $v_{1},\ldots,v_{m}\in V_{N}=\mathbb{F}_q^{~N}$
on which the (right) action of $g$ can be described by a multiplication
from the left by the matrix $\B$. Namely, 
\begin{eqnarray}
\btil_N\left(g\right)&=&\#\left\{ M\in \M\,\middle|\,Mg=\B M\right\} \nonumber\\
&=&\#\left\{ \phi\colon V_m \to V_N
\,\middle|\,\phi \text{ is a linear map, and } \phi g=\B \phi\right\},\label{eq:2nd def of B tilde}
\end{eqnarray}
where the composition of linear transformations in  \eqref{eq:2nd def of B tilde} is from left to right\footnote{Our choice to let $g$ act from the right  is immaterial: $\btil$ could have been equivalently defined as the number of $N\times m$ matrices satisfying $gM=MB$ -- see the proof of Theorem \ref{thm: R^fr and P} in Section \ref{subsec:R^fr and P}.}.
\end{defn}

For example, the trivial and only element of $G_0$ defines the constant-1 function. If $\B=(1)\in G_1$, then
$\btil_N=\fix_N$, where $\fix_N\colon G_N \to\mathbb{Z}_{\ge0}$ is the function
counting elements in $V_N$ which are fixed
by a given matrix in $G_N$. For $\B=(\lambda) \in G_1$,
the function $\btil_N(g)$ gives the size of the eigenspace $V_{\lambda}\le V_{N}$
of $g$. If $\B=I_{m}\in G_m$, then $\btil_N\left(g\right)=\fix\left(g\right)^{m}$,
and if 
\[
\B=\left(\begin{array}{ccccc}
 & 1\\
 &  & 1\\
 &  &  & \ddots\\
 &  &  &  & 1\\
1
\end{array}\right)\in G_m,
\]
then $\btil_N(g)=\fix\left(g^m\right)$. As a final example, if $B = \begin{pmatrix} 
\lambda & 0\\
1 & \lambda 
\end{pmatrix} \in G_2$, then 
\begin{equation}
    \btil_N(g)=|\ker((g-\lambda I)^2)|.\label{eq:first example with size of kernel}
\end{equation}

For every $B\in G_m$, $\btil_N$ is a class function on $G_N$ for all $N$: indeed, for $g,x\in G_N$ and $M\in M_{m\times N}(\mathbb{F}_q)$, we have $Mg=BM$ if and only if $(Mx^{-1})\cdot (xgx^{-1}) = B\cdot (Mx^{-1})$. So $\btil\defi\{\btil_N\}_{N\ge0}$ is an element of (or, more precisely, represents an equivalence class in) $\cfq$.

Define the $\mathbb{C}$-subalgebra 
\[
\R\defi\spn_{\mathbb{C}}\left\{ \btil\,\middle|\,\B\in G_m,m\in\mathbb{Z}_{\ge0}\right\} \subseteq{\cal CF}_q
\]
spanned by all functions $\btil$ with $B\in G_m$ for some $m\in\mathbb{Z}_{\ge0}$\footnote{As mentioned above, the actual elements of $\R$ are the equivalence classes $\left[\btil\right]$ in $\cfq$. For ease of notation, we simply write $\btil$. This will also be true throughout the paper, and in particular when we define the remaining three subalgebras of $\cfq$.}. This is indeed a subalgebra as it is closed under multiplication: 
\begin{equation}
    \btil_1\cdot \btil_2=\widetilde{\begin{pmatrix}\B_1 & 0\\
0 & \B_2
\end{pmatrix}}.\label{eq:multiplication in R}
\end{equation}

Some of the functions in $\R$ are equivalent to "stable irreducible characters" of $G_\bullet$. For example, for every $N\ge2$, the action of $G_N$
on the projective space $\mathbb{P}^{N-1}\left(q\right)$ decomposes
into the trivial representation and an irreducible representation of dimension $\frac{q^N-q}{q-1}$, the
character of which we denote by $\chi^{\mathbb{P}}_N$. Inside $\cfq$, the element $\chi^{\mathbb{P}}=\{\chi^\mathbb{P}_N\}_{N\ge2}$ is equivalent to the following element of $\R$:
\begin{equation}
    \chi^{\mathbb{P}}=-1+\frac{1}{q-1}\sum_{\lambda\in \mathbb{F}_q^{*}}\left(\widetilde{(\lambda)}-1\right)\label{eq:chi^P}
\end{equation}
(here $(\lambda)$ is an element of $G_1$).\\

Our initial motivation behind the current paper was to show that the fact that $\chi^\mathbb{P}$ belongs to $\R$ is not a coincidence but rather the rule: we prove that the $\mathbb{C}$-algebra $\R$ is precisely the algebra of ``stable'' characters, namely, the class functions in $\cfq$ spanned by the characters corresponding to finitely generated $\mathrm{VI}$-modules in the sense of \cite{putman2017representation,gan2018representation} (this is the equality $\R=\I$ in Theorem \ref{thm:main}).
Thus, the algebra $\R$ includes all ``natural families'' of finite
dimensional irreducible representations of $G_\bullet$,
for which one could hope that the theory about word-measures in \cite{ernstwest2022word} may apply. In particular, this is important in light of \cite[Conj.~A.4]{puder2023stable}, which connects a new stable invariant of words in free groups (an invariant in the spirit of stable commutator length -- see [Ibid.~,Def.~A.2]), with the expected value of stable irreducible characters of $G_\bullet$ under word measures.

We will use the natural structure of $\R$ as a filtered algebra\emph{}\footnote{A filtered algebra is an algebra $A$ over some field that has an
increasing sequence $\{ 0\} \subseteq F_0\subseteq F_1\subseteq\cdots$
of subspaces such that $A=\bigcup_{i}F_i$ and $F_m\cdot F_r\subseteq F_{m+r}$. \label{fn:filtered alg}} $\left\{ \R_{\le m}\right\}_{m=0}^\infty$, where
\begin{eqnarray*}
    \R_m & \defi & \spn_{\mathbb{C}}\left\{ \btil\,\middle|\,\B\in G_m\right\}, \text{ and} \\
    \R_{\le m} & \defi & \sum_{j=0}^m\R_j
\end{eqnarray*}
are linear subspaces of $\R$. The filtration $\left\{ \R_{\le m}\right\}_{m=0}^\infty$ yields indeed a filtered algebra, because by \eqref{eq:multiplication in R} we see that $\R_{m}\cdot\R_{r}\subseteq\R_{m+r}$
and $\R_{\le m}\cdot\R_{\le r}\subseteq\R_{\le m+r}$. 

The easy Lemma \ref{lem:conjugates give the same function} below shows that if $B_1,B_2\in G_m$ are conjugate, then $\btil_1=\btil_2$. In Section \ref{sec:immediate corollaries} we show that this is the only source for linear dependencies in $\R_m$, and, in fact, even in $\R$ entirely. We conclude there  that $\left\{ \R_{m}\right\} _{m}$
gives a graded-algebra structure on $\R$, namely, that in addition to the above, $\R=\bigoplus_{m=0}^\infty \R_m$. In Appendix \ref{sec:btil is size of kernel} we also generalize the interpretation in \eqref{eq:first example with size of kernel} to all $\btil$.

%%%%%%%%%%%%%%%%%%%%%%%%%%%%%%%%%%%%%%%%%%%%%%%%%%%%%%%%%%%%%
%   R full rank
%%%%%%%%%%%%%%%%%%%%%%%%%%%%%%%%%%%%%%%%%%%%%%%%%%%%%%%%%%%%%

\subsection*{The subspaces $\R^{\fr}$, $\R_m^\fr$ and $\R_{\le m}^\fr$}

We may define slightly different elements of ${\cal CF}_q$ by restricting
to matrices $M \in \M$ of the maximal possible rank -- \textbf{full rank} $m$. Indeed, for $\B\in G_m$ and $g\in G_N$, denote
\begin{eqnarray}
\btil_N^\fr\left(g\right)&\defi&\#\left\{ M\in \M\,\middle|\,Mg=\B M, ~\rk(M)=m\right\} \nonumber\\
&=&\#\left\{ \phi\colon V_m \hookrightarrow V_N
\,\middle|\,\phi \text{ is an injective linear map, and } \phi g=\B \phi\right\}.\label{eq:def of B tilde full rank}
\end{eqnarray}
and $\btil^\fr\defi\{\btil^\fr_N\}_{N\ge0}$.

The function $\btil^\fr$ was also defined in \cite[p.~1044]{gadish2020dimension} (up to some multiplicative constant depending on $B$). Note that $\btil^\fr_N\equiv0$ for every $N<m$ (this is, of course, irrelevant to the equivalence class of $\btil^\fr$ in $\cfq$).
For $\B\in G_{0}$ the trivial (and only) element, $\btil^{\fr}=\btil$
is the constant function $1$. For $\B=(1) \in G_{1}$,
while $\btil_N=\fix_N$, we have that $\btil^\fr_N=\fix_N-1$ counts the
number of non-zero fixed vectors of an element in $G_{N}$.

In parallel to the above, consider the following subspace of ${\cal CF}_q$:

\[
\R^{\fr}\defi\spn_{\mathbb{C}}\left\{ \btil^{\fr}\,\middle|\,\B\in G_m ,m\in\mathbb{Z}_{\ge0}\right\} \subseteq{\cal CF}_q,
\]    
with the filtrations 
\[
\R_{m}^{\fr}\defi\spn_{\mathbb{C}}\left\{ \btil^{\fr}\,\middle|\,\B\in G_{m}\right\} ~~~~~\mathrm{and}~~~~~\R_{\le m}^{\fr}\defi\sum_{j=0}^{m}\R_{j}^{\fr}.
\]

It is not immediate that $\R^\fr$ is an algebra, namely, closed under products, but it does follow from Gadish's general theory of representation stability and character polynomials \cite[Cor.~3.9]{gadish2017categories}. It also follows immediately from Theorem \ref{thm:main}, which states, in particular, that $\R^\fr=\R$. In Section \ref{sec: Nir} we also prove a conjecture appearing in \cite[p.~979]{balachandran2023product}, which concerns the multiplication in $\R^\fr$.

%%%%%%%%%%%%%%%%%%%%%%%%%%%%%%%%%%%%%%%%%%%%%%%%%%%%%%%%%%%%%
% P
%%%%%%%%%%%%%%%%%%%%%%%%%%%%%%%%%%%%%%%%%%%%%%%%%%%%%%%%%%%%%

\subsection*{The subspaces $\P$, $\P_m$ and $\P_{\le m}$ of parabolic representations}

Although many of the terms here were originally introduced in \cite{zelevinsky1981representations}, we will follow the notations of \cite{gan2018representation} which are more concrete and accessible in our particular case. Given
$m,r\in\mathbb{Z}_{\ge0}$ with $N=m+r$, let $P_{m,r}$ be the (parabolic)
subgroup of $G_{N}$ defined as 
\begin{equation}
P_{m,r}=\left\{ \begin{pmatrix}g_{11} & g_{12}\\
0 & g_{22}
\end{pmatrix}:g_{11}\in G_{m},g_{22}\in G_{r}\right\} .\label{parabolic_subgroup}
\end{equation}
Define the subgroups $G_{m,r}$ and $U_{m,r}$ of $P_{m,r}$ by the
conditions that for any element $p$ of the form \eqref{parabolic_subgroup},
one has: 
\[
\begin{split}p\in G_{m,r} & \iff g_{12}=0,\\
p\in U_{m,r} & \iff g_{11}=I_{m}\textrm{ and }g_{22}=I_{r},
\end{split}
\]
where for any $m\in\mathbb{Z}_{\ge0}$, we write $I_{m}$ for the
identity element of $G_{m}$. The subgroup $G_{m,r}$ is sometimes called a \textit{Levi Factor}. The composition $G_{m,r}\hookrightarrow P_{m,r}\twoheadrightarrow P_{m,r}/U_{m,r}$
is an isomorphism. If $\pi_{1}$ is a representation of $G_{m}$ and
$\pi_{2}$ is a representation of $G_{r}$, the external tensor product representation $\pi_{1}\boxtimes\pi_{2}$ is a representation of $G_{m}\times G_{r}\cong G_{m,r}\cong P_{m,r}/U_{m,r}$.
Denote by $\left(\pi_{1}\boxtimes\pi_{2} \right)^*$ the pullback of $\pi_{1}\boxtimes\pi_{2}$ to a representation of $P_{m,r}.$ By inducing $\left(\pi_{1}\boxtimes\pi_{2} \right)^*$ from
$P_{m,r}$ to $G_{m+r}=G_{N}$ we get a representation of $G_{N}$
which is denoted by $\pi_{1}\circ\pi_{2}$.
By abuse of notation, if $\chi_{1}$ is the character of $\pi_{1}$
and $\chi_{2}$ the character of $\pi_{2}$, we denote by $\left(\chi_{1}\boxtimes\chi_{2} \right)^*$ the character of $\left(\pi_{1}\boxtimes\pi_{2} \right)^*$ and by $\chi_{1}\circ\chi_{2}$
the character of $\pi_{1}\circ\pi_{2}$. (This notation goes back to the classical paper of Green \cite{green1955characters}.)

Denote by $\triv_N$ the trivial character of $G_N$ (by abuse of notation we may think of it both as the trivial representation and as the trivial character). Every character $\chi$ of $G_{m}$ yields a class function $\chi\circ\triv_{N-m}$ on $G_N$ for all $N\ge m$, which is sometimes called a Weil character. We denote the resulting element of $\cfq$ by 
$\chi\circ1$.
For every $m\in\mathbb{Z}_{\ge0}$, consider the subspace 
\[
    \P_m \defi \spn_{\mathbb{C}}\left\{ \chi\circ1\,\middle|\,\chi\in\text{IrrChar}(G_m)\right\} \subseteq{\cal CF}_q,
\]
where $\text{IrrChar}(G_{m})$ denotes the set of irreducible (complex)
characters of $G_{m}$. Define also 
\[
    \P_{\le m}\defi\sum_{j=0}^m \P_j,~~~~\text{and}~~~~
    \P\defi\bigcup_{m\ge0}\P_{\le m}.
\]

For example, ${\cal P}_0$ consists of the constant functions. For any non-trivial $\chi\in\text{IrrChar}(G_1)$ and $N\ge1$, $\chi\circ1$ specializes to the class function $\chi\circ\triv_{N-1}$ of $G_N$, which is a $\frac{q^N-1}{q-1}$-dimensional irreducible character of $G_N$. Distinct non-trivial characters of $G_1$ yield distinct characters of $G_N$. When $\chi=\triv_1$ is the trivial character of $G_1$, the character $\triv_1\circ1$ specializes to the $\frac{q^N-1}{q-1}$-dimensional character of $G_N$ that decomposes, for $N\ge2$, as $\triv_N \oplus \chi^{\mathbb{P}}_N$. So ${\cal P}_{\le1}$ is the $\mathbb{C}$-span of $q$ distinct (what we shall shortly call "stable") irreducible characters, including the trivial one.

%%%%%%%%%%%%%%%%%%%%%%%%%%%%%%%%%%%%%%%%%%%%%%%%%%%%%%%%%%%%%
% I
%%%%%%%%%%%%%%%%%%%%%%%%%%%%%%%%%%%%%%%%%%%%%%%%%%%%%%%%%%%%%

\subsection*{The subspaces $\I$, $\I_m$ and $\I_{\le m}$ spanned by stable irreducible characters}

Here we follow, again, the description in \cite{gan2018representation},
which, in turn, follows \cite{zelevinsky1981representations,springer1984characters}.
An irrep (irreducible representation) $\rho$ of $G_n,$ for $n \in \mathbb{Z}_{\geq 1},$ is
called \textbf{cuspidal} if it does not appear as an irreducible component
of a parabolic representation $\pi_1\circ\pi_2$ for any irreps $\pi_1$ of $G_m$ and $\pi_2$  of $G_{n-m}$, for some $1\le m\le n-1$. In particular, the trivial representation $\triv_1$ of $G_1$ is considered to be cuspidal, while we do not consider $G_0$ to be part of this definition. 
We let $\mathscr{C}_n$ denote the set of all cuspidal irreps of
$G_n$ for $n\ge1$, and $\mathscr{C}\defi\bigsqcup_{n=1}^{\infty}\mathscr{C}_{n}$.
For every $\rho\in\mathscr{C}_{n}$ we define $d(\rho)\defi n$.

The following gives a parametrization of all irreps of $G_{n}$, as
in \cite[p.~49]{gan2018representation}. A \textbf{partition} $\mu=(\mu_1,\mu_2,\ldots,\mu_\ell)$
of $n\in\mathbb{Z}_{\ge0}$ is a finite non-increasing sequence of positive integers with $|\mu|\defi\sum_{i=1}^\ell\mu_{i}=n$. The fact that $|\mu|=n$ is also written as $\mu \vdash n$.
In particular, the only partition of $n = 0$ is the empty one. Let
$\mathscr{P}$ denote the set of all partitions of non-negative integers.
For every function $\vec{\mu}\colon\mathscr{C}\to\mathscr{P}$ with a finite support (so $\vec{\mu}$ maps all but finitely many cuspidal irreps to the empty partition), let
\[
\left\Vert \vec{\mu}\right\Vert \defi\sum_{\rho\in\mathscr{C}}d(\rho)\left|\vec{\mu}\left(\rho\right)\right|.
\]
There is a natural parametrization of the isomorphism classes
of irreps of $G_{n}$ by finitely-supported functions $\vec{\mu}\colon\mathscr{C}\to\mathscr{P}$
such that $\left\Vert \vec{\mu}\right\Vert =n$ \cite[Prop.~9.4]{zelevinsky1981representations}. We denote by $\varphi\left(\vec{\mu}\right)$
the irrep of $G_{n}$ parametrized by $\vec{\mu}$. 

An important role is played by the cuspidal representation $\triv_1$. Suppose that $\vec{\mu}\colon\mathscr{C}\to\mathscr{P}$
satisfies $\vec{\mu}\left(\triv_1\right)=\nu=\left(\nu_{1},\nu_{2},\ldots\right)$,
where, possibly, $\left|\nu\right|=0$ and then we set $\nu_1=0$.
For every integer $N\ge\left\Vert \vec{\mu}\right\Vert +\nu_{1}$,
we define $\vec{\mu}\left[N\right]\colon\mathscr{C}\to\mathscr{P}$
with $\left\Vert \vec{\mu}\left[N\right]\right\Vert =N$ by 
\begin{equation}
\vec{\mu}\left[N\right]\left(\rho\right)=\begin{cases}
\left(N-\left\Vert \vec{\mu}\right\Vert ,\nu_1,\nu_2,\ldots\right) & \textrm{if }\rho=\triv_1,\\
\vec{\mu}\left(\rho\right) & \textrm{if }\rho\neq\triv_1.
\end{cases}
\end{equation}
For example, if $\vec{\mu}$ assigns the empty partition to every cuspidal irrep, then $\vec{\mu}[N]$, defined for all $N\ge0$, assigns the partition $(N)$ to $\triv_1$ and the empty partition to every other cuspidal irrep. In this case, $\varphi(\vec{\mu}[N])=\triv_N$. If $\vec{\mu}$ is supported on $\triv_1$ to which it assigns the partition $(1)$, then $\varphi(\vec{\mu}[N]) = \chi^\mathbb{P}_N$ (these two examples can be inferred from the proof of Proposition \ref{prop:P and I}). 
Clearly, for each function $\vec{\nu}\colon\mathscr{C}\to\mathscr{P}$
with $\left\Vert \vec{\nu}\right\Vert =N<\infty$, there exists a
unique $\vec{\mu}\colon\mathscr{C}\to\mathscr{P}$ such that $\vec{\nu}=\vec{\mu}\left[N\right]$. 

Let $\chi_{\vec{\mu}}$ denote the character of the
representation $\varphi\left(\vec{\mu}\right)$. Following \cite[Def.~1.5]{gan2018representation}, we say that an element of ${\cal CF}_q$ of the form 
\[
    \chimu \defi \{\chi_{\vec{\mu}[N]}\}_{N\ge\left\Vert \vec{\mu}\right\Vert +\nu_{1}}
\]
is a \textbf{stable irreducible character}. For $m\in\mathbb{Z}_{\ge0}$ denote 
\[
\I_m=\spn_{\mathbb{C}}\left\{ \chimu \, \middle| \, \left\Vert \vec{\mu}\right\Vert =m\right\} \subseteq{\cal CF}_q.
\]
For example, $\I_0$ is the span of the trivial character, namely, it is the space of constant functions. Note that for every
$N\ge 2m$, all the elements $\left\{ \varphi\left(\vec{\mu}\left[N\right]\right)\,\middle|\,\left\Vert \vec{\mu}\right\Vert =m\right\} $
are defined and correspond to \emph{distinct} irreps of $G_{N}$, so the set $\left\{ \chimu \,\middle|\,\left\Vert \vec{\mu}\right\Vert =m\right\} $
is a linearly independent set in ${\cal CF}_q$, and 
\begin{equation} \label{dim of Im}
    \dim\left(\I_{m}\right)=\left|\left\{ \vec{\mu}\,\middle|\,\left\Vert \vec{\mu}\right\Vert =m\right\} \right|=\left|\text{IrrChar}(G_m)\right|.
\end{equation}
For a similar reason, for $r\ne m$ we have
$\I_m \cap \I_r=\{0\}$, 
so we may define 
\begin{equation} \label{I as a direct sum}
    \I_{\le m}\defi\bigoplus_{j=0}^{m}{\cal I}_{j}~~~~~ \text{and} ~~~~~ {\cal I}\defi\bigoplus_{m=0}^\infty \I_m = \bigcup_{m=0}^\infty\I_{\le m}.
\end{equation}

We remark that if $\left\Vert \vec{\mu}\right\Vert=m$ then the dimension of  $\varphi(\vec{\mu}[N])$ coincides with $h(q^N)$ where $h\in \mathbb{Q}[x]$ is a polynomial of degree $m$ -- see Corollary \ref{cor:dim of stable irreps}.
In the terminology of \cite{GURALNICK_LARSEN_TIEP_2020}, the irrep $\varphi(\vec{\mu}[N])$ and its character are said to have \textbf{true level} $\left\Vert \vec{\mu}\right\Vert$.

%%%%%%%%%%%%%%%%%%%%%%%%%%%%%%%%%%%%%%%%%%%%%%%%%%%%%%%%%%%%%
%  Main Result
%%%%%%%%%%%%%%%%%%%%%%%%%%%%%%%%%%%%%%%%%%%%%%%%%%%%%%%%%%%%%

\subsection*{The main results}

The following theorem is the main result of the current paper.
\begin{theorem} \label{thm:main} 
The four filtrations
above are all identical, namely, for every $m\in\mathbb{Z}_{\ge0}$,
\[
\R_{\le m}=\R_{\le m}^\fr=\P_{\le m}=\I_{\le m}.
\]
In particular, $$\R=\R^\fr = \P=\I.$$
\end{theorem}

The two middle chains of subspaces, $\R^\fr$ and $\P$, are identical not only on the level of the aggregated, filtered, subspaces, but also in the individual layers, with a nice Fourier-like formula connecting the two bases:
\begin{theorem} \label{thm: R^fr and P}
    For every $m\in\mathbb{Z}_{\ge0}$, 
    \[
        \R^\fr_m = \P_m,
    \]
    with a change of basis given by the following formulas:
    \begin{eqnarray}
        \chi\circ1 &=& \frac1{\left|G_m\right|}\sum_{B\in G_m}\chi(B)\cdot \btil^\fr \label{eq:fourier-like formula 1}\\
        \btil^\fr &=& \sum_{\chi\in \text{IrrChar}(G_m)} \overline{\chi(B)}\cdot (\chi\circ1) \label{eq:fourier-like formula 2}
    \end{eqnarray}
\end{theorem}

We also get interesting and useful formulas which express the generators of $\R_m$ in terms of the generators of $\R_{\le m}^\fr$ and vice versa (Formula \eqref{eq:btil in terms of full rank}), and the generators of $\P_m$ in terms of the generators of $\I_m$ and vice versa \eqref{eq:P in terms of I}. These formulas use some additional notation which  we introduce in Section \ref{sec:proofs}.

The spaces $\I$ and $\P$ have attracted more attention than $\R$ and $\R^\fr$ have, and it is plausible that the equality $\P=\I$ has already appeared in the literature, yet we could not find a reference. 

Theorem \ref{thm:main} gives four different descriptions of the same filtered $\mathbb{C}$-subalgebra of  ${\cal CF}_q$. In view of this result and the connection to stable characters, we call this subalgebra the algebra of \textbf{stable class functions} of $G_\bullet$, and denote it by ${\cal SCF}_q$. Namely, \[{\cal SCF}_q\defi\R=\R^\fr=\P=\I.\]

The equalities from Theorem \ref{thm:main} are useful in many respects. First, some properties of the filtered algebra are obvious in one description but not in others. For example, the linear dimension of every layer is obvious in the case of $\I_m$ -- see \eqref{dim of Im}. This gives us valuable information about the others, such as recognizing precise linear bases in each description -- see Corollary \ref{cor:bases}. 

Second, the multiplication of different basis elements is rather mysterious in three of the constructions: in the case of $\I$ this is basically the notoriously hard problem of decomposing the tensor product of two irreps; in the case of $\R^\fr$, analysing the product is the main goal of the paper \cite{balachandran2023product}. Yet in $\R$ the multiplication is straightforward \eqref{eq:multiplication in R}. The formulas we obtain for change of bases (e.g., Proposition \ref{prop:R and R^fr}), shed light on the multiplication in the other constructions. This allows us to prove a conjecture from [Ibid.] regarding the structure of $\R^\fr$ -- see Section \ref{sec: Nir}.

Finally, as already mentioned above, the theory of word measures in $G_\bullet$ is developed in \cite{ernstwest2022word} using the functions of the form $\btil$, which define $\R$, yet some results and intriguing conjectures about word measures on different families of groups relate to the expected values of stable irreducible characters (see \cite{puder2023stable}). In the case of $G_\bullet$, these are the basis elements of $\I$. The equality $\R=\I$ is an important bridge in this case. See, for instance, Corollary \ref{cor:expected value under uniform} below.

\subsubsection*{Paper organization and notation}
Section \ref{sec:proofs} contains the proofs of Theorems \ref{thm:main} and \ref{thm: R^fr and P}, and provides concrete formulas connecting $\R$ and $\R^\fr$ as well as $\P$ and $\I$. In Section \ref{sec:immediate corollaries} we deduce some corollaries from the main results, some of which we find at least as interesting as the main theorems. Then, Section \ref{sec: S_n} draws analogies between the ring of stable class functions of $G_\bullet$ and the ring of stable class functions on the symmetric groups. In Section \ref{sec: Nir} we further analyse the formulas connecting the basis elements of $\R$ and those of $\R^\fr$ and prove the above-mentioned conjecture from \cite{balachandran2023product} about the structure of $\R^\fr$.
Finally, Appendix \ref{sec:btil is size of kernel} gives an algebraic and combinatorial interpretations of the functions $\btil$, generalizing \eqref{eq:first example with size of kernel}.\\

If $V$ is a vector space then $W\le V$ denotes a subspace $W$. For $m\in \mathbb{Z}_{\ge1}$ we let $[m]\defi\{1,\ldots,m\}$. For $f\in\cfq$, we let $f|_N$ denote its specialization to $G_N$, (which is, of course, well-defined only for large enough values of $N$). Most of our additional running notation was defined in this Section \S\ref{sec:intro}.

\subsubsection*{Acknowledgements}
We thank Inna Entova-Aizenbud, Nir Gadish, Matan Seidel and Noam Ta Shma for beneficial conversations. This work was supported by the European Research Council (ERC) under the European Union’s Horizon 2020 research and innovation programme (grant agreement No 850956), by the Israel Science Foundation, ISF grants 1140/23, as well as National Science Foundation under Grant No. DMS-1926686.

%%%%%%%%%%%%%%%%%%%%%%%%%%%%%%%%%%%%%%%%%%%%%%%%%%%%%%%%%%%%%
%
%  2. Proof of main result
%
%%%%%%%%%%%%%%%%%%%%%%%%%%%%%%%%%%%%%%%%%%%%%%%%%%%%%%%%%%%%%

% Should mention if the proof gives integer coefficients in the base-changes

\section{Proof of the main results} \label{sec:proofs}

In this section we gradually establish the equalities in the statement of Theorem \ref{thm:main}: we first show that $\R_{\le m} =\R_{\le m}^\fr$ in Proposition \ref{prop:R and R^fr}. Then, in Section \ref{subsec:R^fr and P} we show that $\R_m^\fr = \P_m$ and prove Theorem \ref{thm: R^fr and P}. We conclude with the proof that $\P_{\le m}=\I_{\le m}$ in Proposition \ref{prop:P and I}.

We begin with an observation regarding the elements $\btil\in \R_m$ and $\btil^\fr\in \R^\fr_m$ defined with some $B\in G_m$. We already explained above why $\btil_N$ is a class function on $G_N$ for all $N\ge0$: $Mg=BM$ if and only if $(Mx^{-1})xgx^{-1}=B(Mx^{-1})$ for any $g,x\in G_N$. As $M$ has full rank if and only if $Mx^{-1}$ does, this shows that $\btil^\fr\in \R^\fr_m$ is a class function too. But more is true:

\begin{lem}\label{lem:conjugates give the same function}
If $B_1,B_2\in G_m$ are conjugates then $\btil_1=\btil_2$ and $\btil^\fr_1=\btil^\fr_2$.
\end{lem}
\begin{proof}
    Assume that $B_2=xB_1x^{-1}$ for some $x\in G_m$. The lemma follows immediately from the fact that $Mg=B_1M$ if and only if $(xM)g=xB_1x^{-1}(xM)$, and $\rk(M)=\rk(xM)$.
\end{proof}

%%%%%%%%%%%%%%%%%%%%%%%%%%%%%%%%%%%%%%%%%%%%%%%%%%%%%%%%%%%%%
% R and R-full-rank
%%%%%%%%%%%%%%%%%%%%%%%%%%%%%%%%%%%%%%%%%%%%%%%%%%%%%%%%%%%%%

\subsection{$\R$ and $\R^\fr$}

The following proposition expresses $\btil$ as a finite sum of generators of $\R^\fr$. 

\begin{prop}\label{prop:R and R^fr} Let $B\in G_m$ for some $m\in\mathbb{Z}_{\ge0}$. Then there exist $\ell \ge0$ and $C_1,\ldots,C_\ell$ such that for every $j=1,\ldots,\ell$,  $C_j\in G_{k_j}$ for some $k_j<m$, and
\begin{equation}\label{eq:btil in terms of full rank}
    \btil = \btil^\fr + \sum_{j=1}^\ell \tilde{C_j}^\fr.
\end{equation}    
Therefore, for every $m\in\mathbb{Z}_{\ge0}$ we have $\R_{\le m}=\R^\fr_{\le m}$.
\end{prop}

We stress that the matrices $C_j$ are not necessarily distinct: a given matrix $C$ may show up multiple times among $C_1,\ldots,C_\ell$. In fact, in Section \ref{sec: Nir} below, we discuss in detail these multiplicities when $B$ admits a Jordan canonical form.

\begin{proof}
    The last statement follows immediately from \eqref{eq:btil in terms of full rank} by simple induction together with the fact that $\R_0=\R_0^\fr$. So it is enough to prove \eqref{eq:btil in terms of full rank}.

    Fix $B\in G_m$ and $g\in G_N$. The main idea of the proof is to group the different matrices $M\in \text{Mat}_{m\times N}(\mathbb{F}_q)$ satisfying $Mg=BM$ by their left kernel, namely, by the linear dependencies of their rows. Formally, for $M\in \text{Mat}_{m\times N}(\mathbb{F}_q)$, let $\Omega_M\le \mathbb{F}_q^m$ be a subspace of the space of row vectors of length $m$ so that 
    \[
        \omega\in\Omega_M~~\Longleftrightarrow~~\omega M=0.
    \]
    If $Mg=BM$, then for every $\omega\in \Omega_M$ we have $0=\omega Mg = \omega BM$ so $\omega B\in \Omega_M$. Conversely, if $\omega \in \Omega_M$ then $0=\omega Mg^{-1}=\omega B^{-1}M$ so $\omega B^{-1}\in \Omega_M$. Hence, whenver $Mg=BM$, we have $\Omega_MB=\Omega_M$. 
    
    We show below that for every $\Omega\le \mathbb{F}_q^m$ satisfying $\Omega B = \Omega$, the number of matrices $M$ satisfying $Mg=BM$ and $\Omega_M=\Omega$ is precisely $\tilde{C_\Omega}^\fr$ for some $C_\Omega\in G_{m-\dim(\Omega)}$. Obviously, $\dim(\Omega)=0$ if and only if $\Omega=\{0\}$, and in this case $C_\Omega=B$, so this will establish \eqref{eq:btil in terms of full rank}. We stress that the map $\Omega\mapsto C_\Omega$ is generally not one to one.\\

    So now fix $\Omega\le \mathbb{F}_q^m$ with $\Omega B = \Omega$, and let $M\in \text{Mat}_{m\times N}(\fq)$ satisfy $\Omega_M=\Omega$. We will illustrate the argument with the example of $m=4$ and $\Omega =\spn_\mathbb{C}\{(0,1,0,0),(1,0,2,-1)\}$, so the rows of $M$ are of the form $v_1,0,v_2,v_1+2v_2$ for some linearly independent $v_1,v_2$.
    Denote $k=\rk(M)$ and note that $k+\dim(\Omega)=m$. Let $A\subseteq \{1,\ldots,m\}$ be a subset corresponding to the indices of some $k$ linearly independent rows of $M$. Let $S=S(A)\in \text{Mat}_{k\times m}(\fq)$ be a matrix which is obtained from the $m\times m$ identity matrix $I_m$ by keeping only the $k$ rows with indices corresponding to $A$. In our running example, $A$ can be taken to be $\{1,3\}$, in which case 
    \[
        S=\begin{pmatrix}
            1 & 0 & 0 & 0 \\
            0 & 0 & 1 & 0 
        \end{pmatrix}.
    \]
    So $SM$ gives a $k\times N$ matrix of full rank $k$. Let $T\in \text{Mat}_{m\times k}(\fq)$ be a matrix with the property that $TSM=M$, namely, the $i^\text{th}$ row of $T$ expresses the $i^\text{th}$ row of $M$ in terms of the chosen subset of $k$ linearly independent rows. In our example, 
    \[
        T=\begin{pmatrix}
            1 & 0 \\
            0 & 0 \\
            0 & 1 \\
            1 & 2 
        \end{pmatrix}.
    \]
    Note that $ST=I_k$ is the $k\times k$ identity matrix. Also, note that $S$ and $T$ depend only on $A$, which in turn may be chosen according to $\Omega$, but none of $S$, $T$ or $A$ depend on the specific matrix $M$ we consider. In particular, $TSL=L$ holds for every matrix $L$ with $m$ rows and left kernel $\Omega$. Hence there is a bijection 
    \begin{equation} \label{eq:bijection of certain kernel and full}
        \left\{M\in\text{Mat}_{m\times N}(\fq)\,\middle|\, \Omega_M=\Omega \right\}~~~~~ \Longleftrightarrow ~~~~~\left\{L\in\text{Mat}_{k\times N}(\fq)\,\middle|\, \rk(L)=k \right\}
    \end{equation}
    given by $M\mapsto SM$, $L\mapsto TL$.
    
    Denote
    \[
        C=C_\Omega \defi SBT \in \text{Mat}_{k \times k}(\fq).
    \]
    We claim that in the bijection \eqref{eq:bijection of certain kernel and full} above, we have $Mg=BM$ if and only if $Lg=CL$.  Indeed, if $Mg = BM$ then $$Lg=SMg=SBM=SBTSM=CSM=CL.$$ Conversely, if $Lg=CL$ then
    \[
        Mg = TLg = TCL = TSBTL = TSBM = BM,
    \]
    where the last equality follows from the fact that as $\Omega B = \Omega$ we have $\Omega_{BM}=\Omega_M=\Omega$. Hence the number of solutions $Mg=BM$ with $\Omega_M=\Omega$ is precisely $\tilde{C}^\fr(g)$.

\end{proof}

\begin{example} \label{exam:decomposition of B a single Jordan block}
Let $\lambda\in\mathbb{F}_q^*$ and denote by $B_{\lambda,t}\in G_t$ the matrix which is a single Jordan block with eigenvalue $\lambda$, namely, 
\[
    B_{\lambda,t} = \begin{pmatrix}
            \lambda &  & &  & \\
            1&\lambda &  &  & \\
             &1 & \ddots &  & \\
            & & \ddots & \lambda &  \\ 
            & & &  1& \lambda 
        \end{pmatrix} \in G_t.
\]
Then there are exactly $t+1$ subspaces $\Omega\le \mathbb{F}_q^{~t}$ satisfying $\Omega B_{\lambda,t}=\Omega$: these are the subspaces $\Omega_j$ spanned by the first $j$ vectors of the standard basis $\{e_1,e_2,\ldots,e_t\}$, for any $j=0,\ldots,t$. Using the notation from the proof of Proposition \ref{prop:R and R^fr}, we have $C_{\Omega_j}=B_{\lambda,j}$, so
\[
    \widetilde{B_{\lambda,t}}=\widetilde{B_{\lambda,0}}^\fr + \widetilde{B_{\lambda,1}}^\fr + \ldots + \widetilde{B_{\lambda,t}}^\fr.
\]
(The first summand is simply the constant function 1 given by the trivial and only element of $G_0$.)
\end{example}

In Section \ref{sec: Nir} we extend Example \ref{exam:decomposition of B a single Jordan block} to much more general matrices.

%%%%%%%%%%%%%%%%%%%%%%%%%%%%%%%%%%%%%%%%%%%%%%%%%%%%%%%%%%%%%
% R-full-rank and P
%%%%%%%%%%%%%%%%%%%%%%%%%%%%%%%%%%%%%%%%%%%%%%%%%%%%%%%%%%%%%

\subsection{$\R^\fr$ and $\P$} \label{subsec:R^fr and P}
Recall that Theorem \ref{thm:main} states that $\R^\fr_{\le m} = \P_{\le m}$ for all $m$, but as stated in Theorem \ref{thm: R^fr and P}, these two families of subspaces are actually identical on the level of the individual layers, namely, $\R^\fr_m=\P_m$ for all $m$. We now prove this equality and the two concrete formulas \eqref{eq:fourier-like formula 1} and \eqref{eq:fourier-like formula 2} from Theorem \ref{thm: R^fr and P}.

\begin{proof}[Proof of Theorem \ref{thm: R^fr and P}]
Let $\chi \in \text{IrrChar}(G_m)$ and consider $\chi\circ1\in\P_m$. Recall from Section \ref{sec:intro} that $\chi \circ1=\{\chi\circ\triv_{N-m}\}_N$ is induced from the character $(\chi \boxtimes \triv_{N-m})^{\ast}$ of the parabolic subgroup $P_{m,N-m}$. Hence, for every $N\ge m$ and $g \in G_N$,
\begin{equation} \label{eq:induced character}
\chi \circ\triv_{N-m} (g) = \frac{1}{|P_{m,N-m}|}\sum_{\substack{y \in G_N \\ ygy^{-1} \in P_{m,N-m}}}(\chi \boxtimes \triv_{N-m})^{\ast}(ygy^{-1}).
\end{equation}
For any $y\in G_N$, if $ygy^{-1}\in P_{m,N-m}$, then it has the form 
\[
ygy^{-1}=\begin{pmatrix}
B & \ast \\
0 & \ast
\end{pmatrix}
\]
for some $B\in G_m$. Denote by $I_{N,m}$ the $N\times m$ matrix of zeros and ones, with ones exactly on the main diagonal: $I_{N,m}(i,j)=\delta_{ij}$. Then $ygy^{-1}\in P_{m,N-m}$ if and only if $ygy^{-1}I_{N,m}$ has only zeros below the first $m$ rows, and in this case, we have $ygy^{-1}I_{N,m} = I_{N,m}B$. In other words, $ygy^{-1}\in P_{m,N-m}$ if and only if for some (unique) $B\in G_m$ we have 
\begin{equation} \label{eq:eq with x^-1}
    gy^{-1}I_{N,m}=y^{-1}I_{N,m}B,
\end{equation}
and in this case $(\chi \boxtimes \triv_{N-m})^{\ast}(ygy^{-1})=\chi(B)$.

Denote $Q=y^{-1}I_{N,m}$, so \eqref{eq:eq with x^-1} becomes $gQ=QB$. Note that $Q$ is an $N\times m$ matrix of full rank $m$. For any such matrix $Q$, there are exactly 
\[
    c_{N,m}\defi (q^N-q^m)(q^N-q^{m+1})\cdots (q^N-q^{N-1})
\]
different matrices $y\in G_N$ with $y^{-1}I_{N,m}=Q$. Denote by $\mathbbm{1}_{\phi}$ the function with value $1$ when the condition $\phi$ is satisfied and $0$ otherwise. We can now rewrite \eqref{eq:induced character} as
\begin{eqnarray}
    \chi\circ\triv_{N-m}(g) & = & \frac{c_{N,m}}{\left|P_{m,N-m}\right|} \sum_{B\in G_m} \chi(B) \sum_{Q\in \text{Mat}_{N\times m}(\fq) \text{ of rank } m} \mathbbm{1}_{gQ=QB} \nonumber\\
    & = & \frac{1}{\left|G_m\right|} \sum_{B\in G_m} \chi(B) \sum_{Q\in \text{Mat}_{N\times m}(\fq) \text{ of rank } m} \mathbbm{1}_{Q^tg^t=B^tQ^t} \nonumber\\
    & = & \frac{1}{\left|G_m\right|} \sum_{B\in G_m} \chi(B) \cdot \widetilde{B^t}^\fr(g^t). \label{eq:eq with transposes}
\end{eqnarray}

Finally, since a matrix over an arbitrary field is conjugate to its transpose (e.g., \cite[Thm.~66]{kaplansky2003linear}), we have $\widetilde{B^t}^\fr(g^t) = \btil^\fr(g)$ by Lemma \ref{lem:conjugates give the same function} and the fact that $\btil^\fr$ is a class function. Therefore, 
\begin{equation} \label{eq:formula for chi circ 1}
    \chi\circ1 = \frac{1}{\left|G_m\right|} \sum_{B\in G_m} \chi(B) \cdot \btil^\fr, 
\end{equation}
which is precisely \eqref{eq:fourier-like formula 1}. This proves, in particular, that $\chi\circ1 \in \R_m^\fr$, so $\P_m\subseteq\R_m^\fr$.\\

The converse formula \eqref{eq:fourier-like formula 2} follows from the invertibility of the character table of $G_m$. More precisely, let $A\in G_m$ be arbitrary. By \eqref{eq:formula for chi circ 1}, we have
\begin{eqnarray}
    \sum_{\chi \in \text{IrrChar}(G_m)} \overline{\chi(A)} \cdot (\chi\circ1) & = & \frac{1}{\left|G_m\right|} \sum_{\chi \in \text{IrrChar}(G_m)} \overline{\chi(A)} \sum_{B\in G_m} \chi(B) \cdot \btil^\fr \nonumber \\
    & = & \frac{1}{\left|G_m\right|} \sum_{B\in G_m}  \btil^\fr \sum_{\chi \in \text{IrrChar}(G_m)} \overline{\chi(A)} \cdot \chi(B). \label{eq:inner prod of chararacters} 
\end{eqnarray}
By the orthogonality of irreducible characters, the final summation in \eqref{eq:inner prod of chararacters} vanishes unless $B$ and $A$ are conjugates, in which case it is equal to $\frac{|G_m|}{|A^{G_m}|}$, where $A^{G_m}$ is the conjugacy class of $A$ in $G_m$. Recalling again Lemma \ref{lem:conjugates give the same function}, by which $\tilde{A}^\fr=\btil^\fr$ whenever $A$ and $B$ are conjugates, the right hand side of \eqref{eq:inner prod of chararacters} becomes simply $\tilde{A}^\fr$, so we obtain the desired formula \eqref{eq:fourier-like formula 2}. In particular, $\tilde{A}^\fr \in \P_m$, and $\R_m^\fr \subseteq \P_m$.
\end{proof}

%%%%%%%%%%%%%%%%%%%%%%%%%%%%%%%%%%%%%%%%%%%%%%%%%%%%%%%%%%%%%
% P and I
%%%%%%%%%%%%%%%%%%%%%%%%%%%%%%%%%%%%%%%%%%%%%%%%%%%%%%%%%%%%%
\subsection{$\P$ and $\I$}

Our proof that $\P_{\le m}=\I_{\le m}$ uses Pieri's rule for $G_\bullet$. The classical Pieri's rule is a formula for the decomposition to irreps of the induction of certain representations in the symmetric group, and is a special case of the Littlewood-Richardson rule. The analogous rule for $G_\bullet$ was developed independently by \cite{gan2018representation} and \cite{gurevich2020pieri}. In order to present it, we need the following notation from \cite[Notation 2.5]{gan2018representation}.

Recall the notation from Section \ref{sec:intro}. Let $\nu,\mu \in \mathscr{P}$ and $r \in \mathbb{Z}_{\geq 0}.$ We write $\nu  \sim \mu + r$ if the Young diagram of $\nu$ can be obtained by adding $r$ boxes to the Young diagram of $\mu$ with no two boxes added in the same column. Similarly, let $\vec{\nu}, \vec{\mu} : \mathscr{C} \to \mathscr{P}$ and $r \in \mathbb{Z}_{\geq 0}.$ We write $\vec{\nu} \sim \vec{\mu} + r$ if $\vec{\nu}\left(\triv_1 \right) \sim \vec{\mu} \left(\triv_1 \right) + r$ and $\vec{\nu}\left(\rho \right) = \vec{\mu}\left(\rho \right)$ for all $\triv_1 \neq \rho \in \mathscr{C}$. 

\begin{lem}[Pieri's rule for $\gl_N(\mathbb{F}_q)$]\cite[Lem.~2.8]{gan2018representation}, \cite[\S4.5]{gurevich2020pieri} \label{lem:Pieri}
For any $r \in \mathbb{Z}_{\geq 0}$ and a finitely-supported $\vec{\mu} : \mathscr{C} \to \mathscr{P}$,
\begin{equation}\label{eq:Pieri}
    \chi_{\vec{\mu}} \circ \triv_r = \sum_{\vec{\nu} \sim \vec{\mu} + r} \chi_{\vec{\nu}}.
\end{equation}
\end{lem}

\begin{prop} \label{prop:P and I}
For every $m\in\mathbb{Z}_{\ge0}$, $\P_{\le m}=\I_{\le m}$.
\end{prop}

\begin{proof}
For the purpose of this proof, it will be useful to think of $\nu \sim \mu + r$ as being obtained from $\mu$ by subtracting boxes from $\mu$ and then adding an extra row on top, rather than being obtained by adding $r$ boxes to $\mu.$ Formally, let $\mu = \left(\mu_1,\mu_2,\ldots \right)$ and assume that $\nu$ is obtained by adding $r$ boxes to $\mu,$ at most one to each column. Denote by $k \leq \mu_1$ the number of boxes added on to the non-empty columns of $\mu,$ and denote by $J = \left\{j_1,j_2, \ldots, j_k \right\}$ the set of indices of the non-empty columns to which a box was added. So $\nu$ can equivalently be obtained from $\mu$ by subtracting $\mu_1 - k$ boxes from the columns indexed by $\left[\mu_1 \right] \setminus J$ (recall that $[\mu_1]=\{1,\ldots,\mu_1\}$), and then adding a row on top with $r + \mu_1 - k$ boxes.     
\\
\\
For example, if $\mu = \left(5,2,1,1 \right)$ and $\nu = \left(7,4,1,1,1 \right),$ then $\nu  \sim \mu + 5$ and can be obtained from $\mu$ by adding a box to the columns of $\mu$ indexed by $1,3,4,6,7.$ This can be illustrated by the left side of the following diagram:
\begin{figure}[h]
$
\ydiagram[*(red)]
{5+2,2+2,1+0,1+0,1}
*[*(white)]{7,4,1,1,1}
~~~~~~~~~~~~~~~~~~~~~~~~~~\ydiagram[*(red) ]{7}
*[*(blue)]
{7,4+1,1+1,1+0,1+0}
*[*(white)]{7,5,2,1,1}
$
\centering
% \caption{This is the Young diagram of $\nu.$ The boxes in white represent the Young diagram of $\mu$ and the boxes in red are the boxes added in order to obtain the Young diagram of $\nu.$}
\end{figure}

Alternatively, $\nu$ can be obtained from $\mu$ by subtracting a box from columns $2$ and $5$. and adding a row on top with $5 + 5 - 3$ boxes. This is illustrated by the right diagram above.

So in principle, instead of adding $r$ boxes to $\mu$, at most one in each column, we may remove $s$ boxes for arbitrary $s$, at most one from each column, and then add a first row of size $s+r$ as long as we get a valid diagram. One subtlety here is that if $\mu_1>\mu_2$ and we remove the rightmost box of the first row, we still need to add a first row of length at least $\mu_1$, even though adding a shorter first row might create a valid Young diagram. This, however, cannot happen if $r\ge\mu_1-1$.

Recall that in our notation, if $k\ge \vec{\mu}(\triv_1)_1$, then adding a first row of size $k$ to $\vec{\mu}(\triv_1)$ is denoted by $\vec{\mu}[\|\vec{\mu}\| + k]$. This allows us to rewrite Pieri's rule \eqref{eq:Pieri} as follows. For any finitely supported $\vec{\mu} : \mathscr{C} \to \mathscr{P}$ and $r \in \mathbb{Z}_{\geq \vec{\mu}(\triv_1)_1}$,
\[
    \chi_{\vec{\mu}} \circ \triv_r = \sum_{\substack{s \in \mathbb{Z}_{\geq 0} \\ \vec{\nu}\colon\vec{\mu} \sim \vec{\nu} + s }} \chi_{\vec{\nu} \left[\|\mu \| + r\right]}.
\]
Equivalently, if $\| \vec{\mu} \| = m$ and $N\ge m+\vec{\mu}(\triv_1)_1$ we get

\begin{equation}\label{eq:P in terms of I for each N}
    \chi_{\vec{\mu}} \circ \triv_{N-m} = \sum_{\substack{s \in \mathbb{Z}_{\geq 0} \\ \vec{\nu}\colon\vec{\mu} \sim \vec{\nu} + s}}
    \chi_{\vec{\nu}[N]}.
\end{equation}
In our notation for the elements of $\cfq$, this becomes
\begin{equation}\label{eq:P in terms of I}
    \chi_{\vec{\mu}} \circ 1 = \sum_{\substack{s \in \mathbb{Z}_{\geq 0} \\ \vec{\nu}\colon\vec{\mu} \sim \vec{\nu} + s}}
    \chi_{\vec{\nu}[\bullet]}.
\end{equation}
Note that in the right hand side of \eqref{eq:P in terms of I}, every $\vec{\nu}$ satisfies $\|\vec{\nu}\| = \|\vec{\mu}\|-s \le \|\vec{\mu}\|$, so $\chi_{\vec{\nu}[\bullet]} \in \I_{\le m}$. 
This completes the proof that $$\P_{\le m} \subseteq \I_{\le m}.$$ 

We prove the converse by induction on $m$. In the base case $m = 0$, the fact that $\I_0=\P_0$ both consist of the constant functions was already mentioned in Section \ref{sec:intro}. Now assume by the induction hypothesis that $\I_{\le m-1} \subseteq \P_{\le m-1}$. In the right hand side of \eqref{eq:P in terms of I}, the only $\vec{\nu}$ with $\|\vec{\nu}\| = m$ is $\vec{\nu}=\vec{\mu}$. The remaining ones all satisfy $\|\vec{\nu}\| \le m-1$ and so $\chi_{\vec{\nu}[\bullet]}$ belongs to $\I_{m-1} \subseteq \P_{\le m-1} \subseteq \P_{\le m}$. Hence,
\begin{equation} \label{eq:recursion for I in terms of P}
    \chimu = 
    \chi_{\vec{\mu}} \circ 1 - 
    \sum_{\substack{s \in \mathbb{Z}_{\geq 1} \\ \vec{\nu}\colon\vec{\mu} \sim \vec{\nu} + s}}
    \chi_{\vec{\nu}[\bullet]} \in \P_{\le m}.
\end{equation}    
This proves that $\I_m\subseteq\P_m$ and completes the proof of the proposition.
\end{proof}

%%%%%%%%%%%%%%%%%%%%%%%%%%%%%%%%%%%%%%%%%%%%%%%%%%%%%%%%%%%%%
%
%  3. Corollaries and Analogies
%
%%%%%%%%%%%%%%%%%%%%%%%%%%%%%%%%%%%%%%%%%%%%%%%%%%%%%%%%%%%%%

\section{Corollaries of the main result} \label{sec:immediate corollaries}

In this section we describe some basic features of the four families of linear generators of stable class functions of $G_\bullet$ discussed in this paper.
We already mentioned above that the stable irreducible characters $\chi_{\vec{\mu}[\bullet]}$ over all finitely supported $\vec{\mu}\colon\mathscr{C}\to\mathscr{P}$ form a linear basis for $\I$, because if $\vec{\mu_1}\ne\vec{\mu_2}$, then for large enough values of $N$, $\chi_{\vec{\mu_1}[N]}$ and $\chi_{\vec{\mu_2}[N]}$ represent distinct irreducible characters of $G_N$. The following corollary explains how it follows from Theorem \ref{thm:main} that the other three filtrations of ${\cal SCF}_q$ admit, too, natural linear bases. Denote by $\text{ConjClaReps}(G)$ a set of representatives of the conjugacy classes of the group $G$.

\begin{cor} \label{cor:bases}
The space $\scfq$ of stable class functions on $G_\bullet$ admits the following four linear bases:
\begin{enumerate}
    \item $X_m\defi\left\{ \btil \, \middle| \, B \in\text{ConjClaReps}(G_m) \right\}$ is a linear basis of $\R_m$, and $\bigsqcup_{m=0}^\infty X_m$ is a linear basis of ${\cal SCF}_q$.
    
    \item $Y_m\defi\left\{ \btil^\fr \, \middle| \, B \in\text{ConjClaReps}(G_m) \right\}$ is a linear basis of $\R^\fr_m$, and 
    $\bigsqcup_{m=0}^\infty Y_m$ is a linear basis of ${\cal SCF}_q$.

    \item $Z_m\defi\left\{ \chi\circ1 \, \middle| \, \chi\in\text{IrrChar}(G_m)\right\}$ is a linear basis of $\P_m$, and $\bigsqcup_{m=0}^\infty Z_m$ is a linear basis of ${\cal SCF}_q$.  

    \item $W_m\defi\left\{\chimu \, \middle| \, \vec{\mu}\colon\mathscr{C}\to\mathscr{P},~ \left\Vert \vec{\mu}\right\Vert=m\right\}$ is a linear basis of $\I_m$, and $\bigsqcup_{m=0}^\infty W_m$ is a linear basis of ${\cal SCF}_q$.  
\end{enumerate}
\end{cor}

\begin{proof}
    We already explained above the fourth item. It is clear from the definitions of $\R_m$ and $\R_{\le m}$ and from Lemma \ref{lem:conjugates give the same function} that $X_m$ spans $\R_m$ and that its image spans the quotient space $\nicefrac{\R_{\le m}}{\R_{\le m-1}}$. It follows from Theorem \ref{thm:main} and \eqref{I as a direct sum} that 
    \[
        \nicefrac{\R_{\le m}}{\R_{\le m-1}} = \nicefrac{\I_{\le m}}{\I_{\le m-1}} \cong \I_m,
    \]
    and by \eqref{dim of Im}, the dimension of this quotient space is equal to the number of irreducible representations of $G_m$, thus equal to the number of conjugacy classes of $G_m$. Hence the image of $X_m$ does not only span $\nicefrac{\R_{\le m}}{\R_{\le m-1}}$, but is, in fact, a basis of this quotient space. Hence the elements of $X_m$ are linearly independent, they constitute a linear basis of $\R_m$, and $\R_{\le m}=\R_{\le m-1}\oplus \R_m$.

    An identical argument yields the second item of the corollary and a similar one yields the third.
\end{proof}

We already mentioned that $\{\R_{\le m}\}_m$ gives ${\cal SCF}_q$ the structure of a filtered algebra (see Footnote \ref{fn:filtered alg}), and by Theorem \ref{thm:main} this same filtration is given by the other three constructions. Implicit in Corollary \ref{cor:bases} is the following strengthening of this fact:

\begin{cor} \label{cor:R as graded algebra}
    The $\mathbb{C}$-algebra ${\cal SCF}_q$ of stable class functions decomposes as direct sums 
    \[
        {\cal SCF}_q = 
        \bigoplus_{m=0}^\infty \R_m = \bigoplus_{m=0}^\infty \R^\fr_m = \bigoplus_{m=0}^\infty \P_m = \bigoplus_{m=0}^\infty \I_m.
    \]
    In particular, $\{\R_m\}_m$ gives ${\cal SCF}_q$ the structure of a graded algebra. 
\end{cor}

(The filtration $\{\R_m\}_m$ is unique among the four in that $\R_m\cdot \R_\ell \subseteq \R_{m+\ell}$.)\\

The following corollary revolves around the expected value of elements of ${\cal SCF}_q$ under the uniform measure on $G_N$. 
\begin{cor}\label{cor:expected value under uniform}
    For every $f\in{\cal SCF}_q$, the expected value $\mathbb{E}_N[f]$ of $f|_N$ on $G_N$ stabilizes for large enough $N$.
    
    Moreover, for the basis elements from Corollary \ref{cor:bases}, the stabilized average, denoted $\mathbb{E}_\infty[f]$, has the following properties for every $m\in\mathbb{Z}_{\ge0}$:
    \begin{itemize}
        \item For every $B\in G_m$, $\mathbb{E}_\infty[\btil]$ is a positive integer\footnote{In fact, this positive integer has a meaning: the proof of the current corollary together with that of Proposition \ref{prop:R and R^fr}, show that this positive integer is equal to the number of distinct linear subspaces $\Delta\le \mathbb{F}_q^{~m}$ which are invariant under multiplication from the right by $B$.}.
        \item For every $B\in G_m$, $\mathbb{E}_\infty[\btil^\fr]=1$.
        \item For $\chi\in \text{IrrChar}(G_m)$, 
        \[
            \mathbb{E}_\infty[\chi\circ1]=
            \begin{cases}
                1 & \chi=\triv_m\\
                0 & \text{otherwise}.
            \end{cases}.
        \]
        \item For every finitely-supported $\vec{\mu}\colon\mathscr{C}\to\mathscr{P}$,
        $$\mathbb{E}_\infty\left[\chimu \right]=\begin{cases}
            1 & \Vert \vec{\mu} \Vert=0\\
            0 & \text{otherwise}.
        \end{cases}.$$

    \end{itemize}
\end{cor}

\begin{proof}
    The first statement about a general $f\in{\cal SCF}_q$ follows from the elaborated statement about any of the four bases of ${\cal SCF}_q$. We prove the statements about the bases in reverse order. First, whenever $\Vert \vec{\mu} \Vert\ge1$ and $\vec{\mu}[N]$ is defined, the representation $\varphi(\vec{\mu}[N])$ is non-trivial irreducible, and thus its average value on $G_N$ is zero. As for the unique $\vec{\mu}$ with $\Vert \vec{\mu} \Vert=0$, the representation $\varphi(\vec{\mu}[N])$ is the one-dimensional trivial irrep, so its character is the constant function with value $1$.

    The statement about $\mathbb{E}_\infty[\chi\circ1]$ follows from Formula \eqref{eq:P in terms of I}, which expresses $\chi\circ1$ as a sum of distinct elements of the form $\chi_{\vec{\nu}[\bullet]}$. Assume that $\chi=\chi_{\vec{\mu}}$ (necessarily $\Vert \vec{\mu} \Vert=m$).
    The trivial stable irrep appears in that summation if and only if $\vec{\mu}$ is supported solely on the trivial irrep $\triv_1\in \mathscr{C}_1$, and $\vec{\mu}(\triv_1)=(m)$, which precisely means that $\varphi(\vec{\mu})$ is the trivial irrep.

    Now consider some $B\in G_m$. Recall the formula \eqref{eq:fourier-like formula 2} from Theorem \ref{thm: R^fr and P}:
    \[
        \btil^\fr = \sum_{\chi\in \text{IrrChar}(G_m)} \overline{\chi(B)}\cdot (\chi\circ1).
    \]
    By the previous paragraph, the only $\chi\in\text{IrrChar}(G_m)$ for which     $\mathbb{E}_\infty[\chi\circ1]\ne0$ is the trivial character $\triv_m$. As $\triv_m(B)=1$, we deduce that $\mathbb{E}_\infty[\btil^\fr]=1$.    
    Finally, by \eqref{eq:btil in terms of full rank}, $\btil$ is equal to a finite sum of elements of the form $\tilde{C}^\fr$.
\end{proof}

\begin{remark}
    The content of Corollary \ref{cor:expected value under uniform} concerning the elements of $\R^\fr$ was recorded in \cite[Thm.~1.1~and~1.2]{gadish2020dimension} (although there, because of a slightly different normalization of the function $\btil^\fr$, the stable value of the average is not identically 1). The content of Corollary \ref{cor:expected value under uniform} concerning the elements of $\R$ can also be directly deduced from the analysis in \cite[\S2]{ernstwest2022word} if one restricts to the case of $w$ being a single-letter word.
\end{remark}

\begin{cor} \label{cor:inner product}
    For any $f_1,f_2\in{\cal SCF}_q$, the inner product $$\langle f_1,f_2\rangle_N\defi \frac1{|G_N|}\sum_{g\in G_N}f_1|_N(g)\overline{f_2|_N(g)}$$ stabilizes for large enough $N$.
\end{cor}

\begin{proof}
    The inner product $\langle f_1,f_2\rangle_N$ is simply the average value of $f_1|_N\cdot\overline{f_2|_N}$ on $G_N$, and ${\cal SCF}_q$ is closed under multiplication, so the stabilization follows from Corollary $\ref{cor:expected value under uniform}$ as long as we show that $\overline{f_2}\in {\cal SCF}_q$. But the fact that ${\cal SCF}_q$ is closed under complex conjugation follows easily from the observation that $\overline{\chi\circ1}=\overline{\chi}\circ1$.
\end{proof}

The phenomenon recorded in Corollary \ref{cor:expected value under uniform} is a special case of a more general phenomenon revolving around word measures. Every word $w$ in the free group $\F_r$ of rank $r$ induces a word map $w\colon G^r\to G$ on every group $G$, obtained by substituting the letters of $\F_r$ by the $r$ given elements of $G$. For example, if $w=abab^{-2}\in \F_2=\F(a,b)$, then $w\colon G^2\to G$ maps $(g,h)\mapsto ghgh^{-2}$. If $G$ is a finite group, $w$ induces a probability measure on $G$, called the $w$-measure, by pushing forward the uniform measure on $G^r$. In particular, if $w=a$ is a single-letter word, then the $w$-measure is simply the uniform measure.

The study of word measures on $G_\bullet$ is the focus of \cite{ernstwest2022word}. For $f\in {\cal CF}_q$, denote by $\mathbb{E}_w[f|_N]$ the expected value of $f|_N$ under the $w$-measure on $G_N$. 

\begin{cor} \label{cor:word-measures on stable class function}
    Let $1\ne w\in \F_r$, and let $f\in {\cal SCF}_q$ be a stable class function on $G_\bullet$. Then there exists a rational function $h_{w,f}\in \mathbb{C}(x)$, so that for every large enough $N$
    \[
        \mathbb{E}_w[f|_N]=h_{w,f}(q^N).
    \]
    Moreover, the degree of $h_{w,f}$ is non-positive, and so $\lim_{N\to\infty}\mathbb{E}_w[f|_N]$ exists and is finite. Finally, if $w=u^d$ where $d\in\mathbb{Z}_{\ge1}$ and $u\in\F_r$ a non-power, then this limit depends only on $d$ and not on $u$. 
\end{cor}

For example, if $w=a^2b^2c^2\in\F_3=\F(a,b,c)$, $q=3$ and $B=(1)\in G_1$, then $h_{w,\btil}=2+\frac{8(x^2-4x+5)}{(x-1)^2(x-3)^2}$, namely,
\[
    \mathbb{E}_w\left[\btil_N\right] = 2 + \frac{8(3^{2N}-4\cdot{3^N}+5)}{(3^N-1)^2(3^N-3)^2}
\]
(this example is taken from [Ibid., Table 1]). Note that if $w\ne1$ is not a proper power, then Corollary \ref{cor:word-measures on stable class function} yields that $\lim_{N\to\infty}\mathbb{E}_w[f|_N]$ is equal to $\mathbb{E}_\infty[f]$ from Corollary \ref{cor:expected value under uniform}.

\begin{proof}
    The statement of the corollary is true for the elements $\btil$ for arbitrary $B\in G_m$ by  [Ibid., Thm.~1.11~and 1.12]. The general result now follows by Theorem \ref{thm:main}.
\end{proof}

Especially interesting is the fact that $\mathbb{E}_w\left[\chi_{\vec{\mu}[N]}
\right]$ coincides with a rational function in $q^N$: this is important for the validity of the statement of \cite[Conj.~A.4]{puder2023stable}.

We end this section with a corollary about the dimensions of stable irreducible representations.
\begin{cor}\label{cor:dim of stable irreps}
    Let $\vec{\mu}\colon\mathscr{C}\to\mathscr{P}$ be finitely supported with $\Vert \vec{\mu} \Vert=m$. Then there exists a polynomial $h\in \mathbb{Q}[x]$ of degree $m$ so that whenever $\varphi(\vec{\mu}[N])$ is defined,  $$\dim\varphi(\vec{\mu}[N])=h(q^N).$$
\end{cor}
\begin{proof}
    For every $\chi\in\text{IrrChar}(G_m)$, the dimensions of the induced representations $\chi\circ \triv_{N-m}$, $N\ge m$, coincide with a polynomial in $q^N$ of degree $m$ as $$\dim(\chi\circ\triv_{N-m})=\dim(\chi)\cdot[G_N\colon P_{m,N-m}]=\dim(\chi)\frac{(q^N-1)(q^N-q)\cdots(q^N-q^{m-1})}{(q^m-1)(q^m-q)\cdots(q^m-q^{m-1})},$$ (recall that $q$ is not a parameter but rather a fixed value).
    
    We proceed by induction on $m$. 
    For $\vec{\mu}$ with $\Vert \vec{\mu} \Vert=m$, Formula \eqref{eq:P in terms of I for each N} yields that the dimension of $\varphi \left(\vec{\mu}[N] \right)$ is $X-Y$, where $X$ is the dimension of $\varphi \left(\vec{\mu} \right)\circ \triv_{N-m}$, namely a polynomial in $q^N$ of degree $m$, and $Y$ is the sum of dimensions of irreps of the form $\varphi \left(\vec{\nu}[N] \right)$ with $\Vert\vec{\nu}\Vert<m$, namely, by the induction hypothesis, a polynomial of degree $<m$. This proves both the case $m=0$ (where no such $\vec{\nu}$ exist), and the induction step.
\end{proof}
For example, if $\vec{\mu}$ is supported on $\triv_1$ and $\vec{\mu}(\triv_1)=(m)$, it follows from \eqref{eq:P in terms of I for each N} that 
\[
    \chi_{\vec{\mu}[N]} = \triv_m\circ\triv_{N-m} - \triv_{m-1}\circ\triv_{N-m+1}.
\]
So $$\dim(\varphi(\vec{\mu}[N]))= \frac{(q^N-1)(q^N-q)\cdots(q^N-q^{m-2})}{(q^m-1)(q^m-q)\cdots(q^m-q^{m-2})}\cdot \frac{q^N-q^m}{q^m-q^{m-1}}.$$

\begin{remark}
    Corollary \ref{cor:dim of stable irreps} can also be deduced from the $q$-analogue of the hook formula -- see \cite[\S11.10]{zelevinsky1981representations}. The fact that the dimension of $\varphi(\vec{\mu}[N])$ agrees with a polynomial in $q^N$ also appears in \cite[Thm.~1.7]{gan2018representation}, but the degree is only shown to be at most $\Vert\vec{\mu}\Vert$ [Ibid., Rem.~5.3].
\end{remark}

%%%%%%%%%%%%%%%%%%%%%%%%%%%%%%%%%%%%%%%%%%%%%%%%%%%%%%%%%%%%%
% S_n
%%%%%%%%%%%%%%%%%%%%%%%%%%%%%%%%%%%%%%%%%%%%%%%%%%%%%%%%%%%%%
\section{Analogies with the symmetric group} \label{sec: S_n}

The algebra of stable class function over $G_\bullet=\gl_\bullet(\fq)$ and its four graded bases considered in the current paper all have analogs over the symmetric group, which is often considered to be the matrix group over the "field of one element". We find it illuminating to briefly describe these analogs and draw the similarities between the two cases. 

Denote by $S_N$ the symmetric group on $N$ elements, and let $S_\bullet \defi \bigsqcup_{N=0}^\infty S_N$.  The algebra of stable class functions of $S_\bullet$, which we denote here by $\scf$, is mentioned in different contexts in the literature. It is hinted to in classical texts such as \cite[Example I.7.14]{macdonald1998symmetric}, and studied more systematically, for example, in \cite{church2015fi}. 

\paragraph{The analog of $\R$}
Recall that for $B\in G_m$, we defined $\btil_N(g)$ as the number of linear maps $\phi\colon V_m\to V_N$ such that $\phi g =B\phi$. The analog of this counting function is the following: given permutations $\sigma\in S_m$ and $g\in S_N$, denote by 
\[
    \tilde{\sigma}_N(g)\defi \#\left\{ \phi\colon [m]\to[N]~\middle|~g\circ\phi=\phi\circ\sigma  \right\},
\]
and $\tilde{\sigma}\defi\{\tilde{\sigma}_N\}_N\in\scf$.

It is not hard to see that if $\sigma$ has cycle structure $\mu=(\mu_1,\ldots,\mu_\ell) \vdash m$, then $$\tilde{\sigma}_N(g)=\fix(g^{\mu_1})\cdots \fix(g^{\mu_\ell}),$$ where here $\fix(g)$ is the number of fixed points of the permutation $g$. Note that these are, in fact, power-sum symmetric functions in the eigenvalues of the permutation matrix corresponding to $g$.

Not surprisingly, this type of stable functions over $S_\bullet$ is the most natural one to approach when studying word measures on $S_\bullet$, and this is indeed the main object of study in \cite{hanany2023word}. 

\paragraph{The analog of $\R^\fr$}
Recall that for $B\in G_m$, $\btil^\fr_N(g)$ is the number of \textit{injective} linear maps $\phi\colon V_m\hookrightarrow V_N$ such that $\phi g =B\phi$. Analogously, for any $\sigma\in S_m$ and $g\in S_N$, denote by 
\[
    \tilde{\sigma}^\fr_N(g)\defi \#\left\{ \phi\colon [m]\hookrightarrow[N]~\middle|~g\circ\phi=\phi\circ\sigma \right\},
\]
and $\tilde{\sigma}^\fr\defi\{\tilde{\sigma}^\fr_N\}_N\in\scf$.

For $\ell\in\mathbb{Z}_{\ge1}$ denote by $a_\ell\in\scf$ the function that returns the number of $\ell$-cycles in a given permutation. The functions $\tilde{\sigma}^\fr$ are a certain type of cycle-statistics functions on $g$:
\begin{claim}
    If $\sigma \in S_m$ has exactly $r_\ell$ $\ell$-cycles for $\ell=1,\ldots,m$, then 
    \begin{equation} \label{eq:tilde sigma}
        \tilde{\sigma}^\fr=\prod_{\ell=1}^m \ell^{r_\ell}(a_\ell)_{r_\ell},
    \end{equation} where $(k)_r\defi k(k-1)\cdots(k-r+1)$ is the falling factorial.    
\end{claim}
For example, if $\sigma=(12)(34)\in S_7$, then $\tilde{\sigma}^\fr=4a_1(a_1-1)(a_1-2)a_2(a_2-1)$.

\begin{proof}
    Given $\sigma\in S_m$ and $g\in S_N$, we ought to count embeddings $\phi\colon [m]\hookrightarrow[N]$ with $g\circ\phi=\phi\circ \sigma$. Let $x\in[m]$ belong to an $\ell$-cycle of $\sigma$. Then the equation $g\circ\phi=\phi\circ \sigma$ means that $\phi(x),\phi(\sigma(x)),\ldots,\phi(\sigma^{\ell-1}(x))$ must form an $\ell$-cycle of $g$ (in that given cyclic order). In other words, the embedding $\phi$ must map every $\ell$-cycle of $\sigma$ into an $\ell$-cycle of $g$. The first $\ell$-cycle of $\sigma$ can be embedded in one of $\ell\cdot a_\ell|_N(g)$ options: there are $\ell$ possible embedding in each of the $a_\ell|_N(g)$ $\ell$-cycles of $g$. The second $\ell$-cycle of $\sigma$ has now $\ell(a_\ell|_N(g)-1)$ possible embeddings, as we cannot pick the same cycle of $g$ twice. Continuing this way gives \eqref{eq:tilde sigma}.
\end{proof}

Let $g\in S_N$ be uniformly random. It is well known that as $N\to\infty$, $a_\ell|_N(g)$ converges to a Poisson random variable with parameter $\frac1\ell$, and, in fact, each moment stabilizes for large enough $N$ (see, e.g., \cite[\S5]{diaconis1994eigenvalues} for the latter fact and for those in the next paragraph). If $Z_\lambda$ is a Poisson random variable with parameter $\lambda$, then $\mathbb{E}[(Z_\lambda)_r]=\lambda^r$. Thus, $\mathbb{E}[\ell^{r_\ell}(a_\ell|_N(g))_{r_\ell}]$ stabilizes to the value 1.

Moreover, the cycle-statistics $a_1|_N(g),a_2|_N(g),a_3|_N(g),\ldots$ are asymptotically independent, and their joint moments stabilize to joint moments of independent Poisson variables. Thus, the expectation of \eqref{eq:tilde sigma} over uniformly random $g\in S_N$ stabilizes to the value 1. This is analogous to the fact from Corollary \ref{cor:expected value under uniform} that $\mathbb{E}_\infty[\btil^\fr]=1$ for all $B\in G_m$.

\paragraph{The analog of $\P$}
The analog of the generators $\chi\circ1$ of $\P$ is functions of the form $\xi\circ1$ where $\xi\in\textrm{IrrChar}(S_m)$. For every $N\ge m$, the specialization of this function to $S_N$ is defined by
\[
    \xi\circ1|_N\defi\xi\circ\triv_{N-m}\defi\textrm{Ind}_{S_m\times S_{N-m}}^{S_N}\left(\xi\boxtimes\triv_{N-m}\right),
\]
where $\triv_{N-m}$ is the trivial character of $S_{N-m}$. 

\paragraph{The analog of $\I$}
Stable irreducible characters of $S_\bullet$ are very similar to their counterparts in the world of $G_\bullet$. Irreducible characters of $S_N$ are in bijection with partitions of $N$. And for every partition $\mu\vdash m$, one gets a partition $\mu[N]$ for every $N\ge m+\mu_1$ by adding to $\mu$ a first row of size $N-m$. In fact, this is completely parallel to the definition of $\chimu\in\scfq$ from Section \ref{sec:intro} if we follow the definition of cuspidal irreps and notice that the only cuspidal irrep in $S_\bullet$ is the trivial irrep of $S_1$.

\paragraph{The analogous results}
All the theorems of this paper are true, and usually known, when adapted to $S_\bullet$. Our main theorem -- Theorem \ref{thm:main} -- is true for $S_N$, with different parts appearing in different papers. The fact that all four sets of functions are bases for the same algebra $\scf$, is recorded, for example, in \cite[Thm.~3.3.4]{church2015fi} ($\R^\fr$ and $\I$), in \cite[Appendix B]{hanany2023word} ($\R$, $\R^\fr$ and $\I$), and in \cite[Prop.~4.23]{puder2023stable} ($\P$ and $\I$).\\

Moreover, the formulas connecting $\R^\fr$ and $\P$ in Theorem \ref{thm: R^fr and P} hold in $S_\bullet$ as well, and with basically the same proof. Because we could not find a reference, we add a proof here.

\begin{prop} \label{prop:Fourier-like formulas for Sn}
    For every $\xi\in \textrm{IrrChar}(S_m)$ 
    \begin{equation} \label{eq:P in terms of R^fr S_n}
        \xi\circ1 = \frac1{m!} \sum_{\sigma\in S_m} \xi(\sigma)\cdot \tilde{\sigma}^\fr,
    \end{equation} \label{eq:R^fr in terms of P S_n}
    and for every $\sigma\in S_m$,
    \begin{equation} \label{eq:S_n Rfr in terms of P}
        \tilde{\sigma}^\fr = \sum_{\xi\in\textrm{IrrChar}(S_m)}\xi(\sigma)\cdot(\xi\circ1).
    \end{equation}
\end{prop}
(In contrast to \eqref{eq:fourier-like formula 2}, $\xi(\sigma)$ does not need to be conjugated in \eqref{eq:S_n Rfr in terms of P}, as the irreducible characters of $S_m$ are all real valued.)

\begin{proof}
    We follow the outline of the proof of Theorem \ref{thm: R^fr and P} in Section \ref{subsec:R^fr and P}. For $\xi\in\textrm{IrrChar}(S_m)$ and every $g\in S_N$,
    \begin{equation} \label{eq:induced character S_N}
        \xi \circ\triv_{N-m} (g) = \frac{1}{|S_m\times S_{N-m}|}\sum_{\substack{y \in S_N \\ ygy^{-1}=(\sigma,1) \in S_m\times S_{N-m}}}\xi(\sigma).
    \end{equation}
    By abuse of notation, think of $\sigma\in S_m$ also as an element of $S_N$ via the inclusion of $S_m$ in $S_N$ as the stabilizer of $\{m+1,\ldots,N\}$. The equation $ygy^{-1}=\sigma$ is equivalent to $gy^{-1}=y^{-1}\sigma$. As $y^{-1}\in S_N$ is a permutation, it induces an injective map $\phi\colon[m]\hookrightarrow[N]$. The number of $y\in S_N$ so that $y^{-1}$ induces a given $\phi\colon[m]\hookrightarrow[N]$ is precisely $(N-m)!$, so the total number of times we get a summand of $\xi(\sigma)$ in \eqref{eq:induced character S_N} is precisely $(N-m)!\cdot \tilde{\sigma}^\fr_N(g)$. This yields \eqref{eq:P in terms of R^fr S_n}.\\

    Now let $\tau\in S_m$ be arbitrary. By \eqref{eq:P in terms of R^fr S_n}, we have
    \begin{eqnarray}
        \sum_{\xi \in \text{IrrChar}(S_m)} \overline{\xi(\tau)} \cdot (\xi\circ1) & = & \frac{1}{m!} \sum_{\xi \in \text{IrrChar}(S_m)} \overline{\xi(\tau)} \sum_{\sigma\in S_m} \xi(\sigma) \cdot \tilde{\sigma}^\fr \nonumber \\
    & = & \frac{1}{m!} \sum_{\sigma\in S_m}  \tilde{\sigma}^\fr \sum_{\xi \in \text{IrrChar}(S_m)} \overline{\xi(\tau)} \cdot \xi(\sigma) \label{eq:inner prod of chararacters S_n} 
\end{eqnarray}
By the orthogonality of irreducible characters, the final summation in \eqref{eq:inner prod of chararacters S_n} vanishes unless $\sigma$ and $\tau$ are conjugates, in which case it is equal to $\frac{|S_m|}{|\tau^{S_m}|}$, where $\tau^{S_m}$ is the conjugacy class of $\tau$ in $S_m$. As $\tilde{\sigma}^\fr=\tilde{\tau}^\fr$ whenever $\sigma$ and $\tau$ are conjugates, the right hand side of \eqref{eq:inner prod of chararacters S_n} becomes simply $\tilde{\tau}^\fr$, so we obtain the desired formula \eqref{eq:S_n Rfr in terms of P}.

\end{proof}

All the remaining corollaries we spelled out in Section \ref{sec:immediate corollaries} hold in the case of $S_\bullet$ and $\scf$ as well. The analogs of Corollaries \ref{cor:bases}, \ref{cor:R as graded algebra}, \ref{cor:expected value under uniform} and \ref{cor:inner product} hold with the same proofs (as mentioned above regarding the expected value of $\tilde{\sigma}^\fr$, most of these results are already known). The $S_\bullet$-analog of Corollary \ref{cor:word-measures on stable class function} states that for every $w\in\F_r$ and every $f\in\scf$ there exists a rational function $h_{w,f}\in \mathbb{C}(x)$ so that for every large enough $N$, $\mathbb{E}_w[f|_N]=h_{w,f}(N)$. This fact is mentioned in the very beginning of \cite[\S1.1]{hanany2023word} and follows from the analysis in that paper.  Finally, the fact analogous to Corollary \ref{cor:dim of stable irreps}, that given $\mu\vdash m$, the dimension of the irreducible representation of $S_N$ corresponding to $\mu[N]$ coincides with a fixed polynomial in $N$ of degree $m$ is well-known and follows easily from the hook formula.

%%%%%%%%%%%%%%%%%%%%%%%%%%%%%%%%%%%%%%%%%%%%%%%%%%%%%%%%%%%%%
% Proof of Nir's conjecture
%%%%%%%%%%%%%%%%%%%%%%%%%%%%%%%%%%%%%%%%%%%%%%%%%%%%%%%%%%%%%
\section{A proof of a conjecture about the structure of $\R^\fr$} \label{sec: Nir}

In the current section we provide further information about the formula \eqref{eq:btil in terms of full rank} expressing a basis element of $\R$ as a combination of basis elements from $\R^\fr$, and prove a conjecture from \cite{balachandran2023product}. This conjecture is the content of the following Theorem \ref{thm:Nir}. Let $\lambda\in\fq^*$ be a non-zero scalar and $\mu=(\mu_1,\ldots,\mu_k)\in\mathscr{P}$ a partition. Denote by $B_{\lambda,\mu}\in G_{|\mu|}$ the matrix with $k$ Jordan blocks of sizes $\mu_1,\ldots,\mu_k$, each with eigenvalue $\lambda$. If the partition has a single block, so $\mu=(t)$ for some $t\in \mathbb{Z}_{\ge0}$, denote $B_{\lambda,t}\defi B_{\lambda,(t)}$, as in Example \ref{exam:decomposition of B a single Jordan block}. Note that unlike in other sections of the current paper, in the following result $q$ is not fixed.

\begin{theorem} \label{thm:Nir}
The following conjecture from \cite[p.~979]{balachandran2023product} is true: for every partition $\mu=(\mu_1,\ldots,\mu_k)$ there is a degree-$k$ polynomial $p_\mu\in \mathbb{Z}[q][x_1,\ldots,x_{\mu_1}]$ on $\mu_1$ variables and coefficients which are polynomials over\footnote{The original conjecture in \cite{balachandran2023product} was stated with a slightly different definition of $\btil^\fr$ (different by some multiplicative constant depending on $B$), and the coefficients of $p_\mu$ were, accordingly, rational functions in $q$ rather than mere polynomials.} $\mathbb{Z}$ in $q$, so that for every prime power $q$ and every $\lambda\in\fq^*$, 
\[
    \widetilde{B_{\lambda,\mu}}^\fr = p_\mu(q)\big(\widetilde{B_{\lambda,1}}^\fr,\ldots,\widetilde{B_{\lambda,\mu_1}}^\fr\big).
\]
\end{theorem}

Note, in particular, that the polynomial $p_\mu$ does not depend on $\lambda$, and its dependence on the size of the field, $q$, is only at the coefficients, which are given by polynomials in $q$. Part of the motivation behind this conjecture comes from it being a $q$-analog of a result about the symmetric groups (see \cite[Thm.~3.3.4]{church2015fi}), which we already mentioned in Section \ref{sec: S_n} above. Theorem \ref{thm:Nir} is a more precise and elaborate $q$-analog than the result stated in Proposition \ref{prop:R and R^fr}. For further motivation, consult \cite{balachandran2023product}.

We prove this conjecture using the connection we established between $\R$ and $\R^\fr$, while also shedding more light on this connection. Recall from Formula \eqref{eq:btil in terms of full rank} that for every $B\in G_m$, 
\begin{equation}\label{eq:again btil in terms of full rank}
    \btil=\btil^\fr + \sum_{j=1}^\ell \tilde{C_j}^\fr
\end{equation} 
with $C_j\in G_\bullet$ a matrix of size strictly smaller than that of $B$ for every $j$. Moreover, by Lemma \ref{lem:conjugates give the same function} and Corollary \ref{cor:bases}, the multiset of $C_j$'s in \eqref{eq:again btil in terms of full rank} is unique up to conjugation of its elements. 
The following proposition characterizes the possible $C_j$'s in \eqref{eq:again btil in terms of full rank}. 

\begin{prop} \label{prop:C as a block in B}
    Let $B\in G_m$ for some $m\in \mathbb{Z}_{\ge0}$ and let $C=C_j$ from \eqref{eq:again btil in terms of full rank} for some $j=1,\ldots,\ell$. Then up to conjugation, $B$ is of the form 
    \begin{equation} \label{eq:constraint on C}
        \begin{pmatrix}
            C      & * \\
            0      & *
        \end{pmatrix}.
    \end{equation}
\end{prop}

\begin{proof}
    Recall from the proof of Proposition \ref{prop:R and R^fr} that $C=C_\Omega$ for some $\Omega\le \fq^m$ with $\Omega B=\Omega$. Denote $k=m-\dim(\Omega)$. Recall that we chose some subset $A\subseteq[m]$ of $k$ indices of independent rows in some (and, in fact, every) matrix $M$ with left kernel $\Omega$, and that according to $A$ we defined two matrices $S$ and $T$ of dimensions $k\times m$ and $m\times k$, respectively, so that $C=SBT$, $ST=I_k$ and $TSM=M$ for every matrix $M$ with $m$ rows and left kernel $\Omega$.

    We extend $T$ to an $m\times m$ matrix $\hat{T}\in G_m$ by spreading the columns of $T$ to be precisely in the indices $A$, and for every $i\notin A$, we let the $i^\text{th}$ column be the $i^\text{th}$ standard unit vector. We denote by $\hat{S}\defi \hat{T}^{-1}$ the inverse of $\hat{T}$.    
    For example, continuing the running example from the proof of Proposition \ref{prop:R and R^fr}, we have 
    \[
        \hat{T}=\begin{pmatrix}
            1 & 0 & 0 & 0 \\
            0 & 1 & 0 & 0 \\
            0 & 0 & 1 & 0 \\
            1 & 0 & 2 & 1
        \end{pmatrix}
        ~~~\textrm{and}~~~~
        \hat{S}=\begin{pmatrix}
            1 & 0 & 0 & 0 \\
            0 & 1 & 0 & 0 \\
            0 & 0 & 1 & 0 \\
            -1 & 0 & -2 & 1
        \end{pmatrix}.
    \]
    The facts that $\hat{T}$ is indeed non-singular and that the rows indexed by $A$ in $\hat{S}$ give rise to the matrix $S$, are both easily verified. So $S$ is obtained from $\hat{S}$ by omitting the rows of indices outside $A$, and $T$ obtained from $\hat{T}$ by omitting the columns with the same property.

    Now notice that by construction, the left kernel of $T$ is $\Omega$, and because $T$ is of rank $k=m-\dim(\Omega)$, its columns span every column of a matrix with left kernel $\Omega$. As the $i^\text{th}$ row of $\hat{S}T$ is a zero row for every $i\notin A$, we get that the same is true for any matrix $M$ with left kernel $\Omega$: the rows indexed by $i\notin A$ are all zero in $\hat{S}M$ as well. 

    Finally, note that as $\Omega B=\Omega$, the left kernel of $BT$ is also $\Omega$. So in $\hat{S}BT$ the rows indexed by $[m]\setminus A$ are all zeros. Note also that $\hat{S}BT$ is obtained from $\hat{S}B\hat{T}$ by omitting the columns indexed by $[m]\setminus A$. We conclude that in $\hat{S}B\hat{T}$, which is a conjugate of $B$, we have:
    \begin{itemize}
        \item the square submatrix indexed by $A\times A$ is $C=SBT$, and
        \item the submatrix indexed by $([m]\setminus A)\times A$ is a zero matrix. 
    \end{itemize}
    This completes the proof of the proposition.
\end{proof}

\begin{cor} \label{cor:C is defined by lambda and partition nu}
    Let $\lambda\in\fq^*$ and $\mu\in \mathscr{P}$. Then in the decomposition $\widetilde{B_{\lambda,\mu}}= \widetilde{B_{\lambda,\mu}}^\fr     + \sum_j\widetilde{C_j}^\fr$, each $C_j$ is of the form $C_j=B_{\lambda,\nu}$ for some $\nu\in \mathscr{P}$. 
\end{cor}
\begin{proof}
    This is immediate from Proposition \ref{prop:C as a block in B}, which shows, in particular, that the characteristic polynomial of $C_j$ divides that of $B$.
\end{proof}

Our next aim is to more accurately  understand the partitions $\nu$ for which $\widetilde{B_{\lambda,\nu}}^\fr$ appears in the decomposition of $\widetilde{B_{\lambda,\mu}}$. We will use the following notations. 
\begin{itemize}
    \item For any partition $\mu$, we denote by $\mu'$ the conjugate partition (corresponding to the transposed Young diagram), so that $\mu'_j$ is the length of the $j^{\text{th}}$ column of $\mu$.
    
    \item For any matrix $E\in\text{Mat}_{m\times m}(\fq)$ (or any linear transformation $T:V_m\to V_m$) and any $S\subseteq V_m=\fq^m$, denote by $(S)_E$ ($(S)_T$, respectively) the smallest subspace of $V_m$ containing $S$ and closed under multiplication from the right by $E$ (under $T$, respectively). Equivalently, this is the subspace spanned by $v,vE,vE^2,\ldots,vE^{m-1}$ (or $v,T(v),\ldots,T^{m-1}(v)$). When $S=\{v_1,\ldots,v_r\}$ we may simply write $(v_1,\ldots,v_r)_E$ (or $(v_1,\ldots,v_r)_T$). \label{def of (v)_B}

\end{itemize}

The following lemma recovers the partition $\mu\vdash m$ from the manner in which $B_{\lambda,\mu}$ acts on the vector space $V_m=\fq^m$. We omit its trivial proof.

\begin{lem} \label{lem:extracting mu from action}
    Let $B=B_{\lambda,\mu}$ for some partition $\mu\vdash m$. Then for all $j$,
    \[
        \mu'_j = \rk \left((B-\lambda)^{j-1}\right) - \rk \left((B-\lambda)^j\right).
    \]
\end{lem}
Below we will approach the rank of the matrix $(B-\lambda)^j$ through the dimension of the image of the linear action of $(B-\lambda)^j$ on $V_m$.

Returning to the proof of Theorem \ref{thm:Nir}, recall that by Corollary \ref{cor:C is defined by lambda and partition nu}, 
\begin{equation} \label{eq:introducing cal N}
    \widetilde{B_{\lambda,\mu}} = \widetilde{B_{\lambda,\mu}}^\fr + \sum_{\nu\in \cal{N}} \widetilde{B_{\lambda,\nu}}^\fr
\end{equation} 
where $\cal{N}=\cal{N_{\lambda,\mu}}$ is a multiset of partitions which depends on $\mu$ and $\lambda$. 

\begin{cor} \label{cor:nu contained in mu}
    For each $\nu\in\cal{N}$, the Young diagram corresponding to $\nu$ is strictly "contained" in that of $\mu$: for every $i$, $\nu_i\le \mu_i$ (possibly adding trailing zeros), and $|\nu|<|\mu|$.
\end{cor}

\begin{proof}
First, $|\nu|< |\mu|$ by the fact from Proposition \ref{prop:R and R^fr} that for every $\nu\in{\cal N}$, $C\defi B_{\lambda,\nu}$ is a matrix of dimension strictly smaller than that of $B\defi B_{\lambda,\mu}$. The claim is now equivalent to that $\nu'_j\le \mu'_j$ for every $j$. By Proposition \ref{prop:C as a block in B}, the action of $C$ on $\fq^{|\nu|}$ is the same as the action of $B$ on some $|\nu|$-dimensional invariant subspace $W$ of $V_m=\fq^{|\mu|}$. By Lemma \ref{lem:extracting mu from action}, $\nu'_j$ is equal to the maximal number $r$ of vectors $w_1,\ldots,w_r\in W$ such that $\big\{w_i(C-\lambda)^{j-1}\big\}_{i=1}^r$ are linearly independent yet $w_i(C-\lambda)^j=0$ for all $i$. Using this interpretation of $\nu'_j$ and recalling that $W$ is a $B$-invariant subspace of $V_m$, it is clear that $\mu'_j$ is at least as large as $\nu'_j$.\end{proof}

\begin{lem} \label{lem:poly degree k}
    Given $\lambda\in \fq^*$ and a partition $\mu=(\mu_1,\ldots,\mu_k)\in\mathscr{P}$, there is a degree-$k$ polynomial $p_{\lambda,\mu}$ with integer coefficients in $\mu_1$ variables so that 
    \[
        \widetilde{B_{\lambda,\mu}}^\fr = p_{\lambda,\mu}\big(\widetilde{B_{\lambda,1}}^\fr,\ldots,\widetilde{B_{\lambda,\mu_1}}^\fr\big).
    \]
\end{lem}

\begin{proof}
    We prove, moreover, that every monomial $\widetilde{B_{\lambda,t_1}}^\fr\cdots \widetilde{B_{\lambda,t_j}}^\fr$ in the polynomial $p_{\lambda,\mu}$ satisfies $t_1+\ldots+t_j\le |\mu|$.
    We prove both claims simultaneously by induction on $|\mu|$. When $\mu=\emptyset$ is the empty partition, $p_{\lambda,\emptyset}=1$ and both claims hold.

    Now fix $\mu=(\mu_1,\ldots,\mu_k)$ with $|\mu|\ge1$. By Proposition \ref{prop:R and R^fr} and Corollaries \ref{cor:C is defined by lambda and partition nu} and \ref{cor:nu contained in mu}, 
    \begin{equation} \label{eq:expression for B_lambda,mu^fr}
        \widetilde{B_{\lambda,\mu}}^\fr = \widetilde{B_{\lambda,\mu}} - \sum_{\nu\in\mathcal{N}} \widetilde{B_{\lambda,\nu}}^\fr,
    \end{equation}
    where $\mathcal{N}=\mathcal{N}_{\lambda,\mu}$ and $\nu$ is strictly contained in $\mu$ for every $\nu\in\mathcal{N}$. By \eqref{eq:multiplication in R} and Example \ref{exam:decomposition of B a single Jordan block}, 
    \begin{equation} \label{eq:poly for B in R}   \widetilde{B_{\lambda,\mu}}=\widetilde{B_{\lambda,\mu_1}}\cdots \widetilde{B_{\lambda,\mu_k}} = \big(\widetilde{B_{\lambda,0}}^\fr+\ldots+\widetilde{B_{\lambda,\mu_1}}^\fr\big) \cdots \big(\widetilde{B_{\lambda,0}}^\fr+\ldots+\widetilde{B_{\lambda,\mu_k}}^\fr\big),
    \end{equation}
    so $\widetilde{B_{\lambda,\mu}}$ is a degree-$k$ polynomial with integer coefficients in\footnote{Recall that $\widetilde{B_{\lambda,0}}^\fr=1$.} $\widetilde{B_{\lambda,1}}^\fr,\ldots,\widetilde{B_{\lambda,\mu_1}}^\fr$ with the sum $t_1+\ldots+t_j$ in every monomial being at most $|\mu|$ and a single monomial --- $\widetilde{B_{\lambda,\mu_1}}^\fr\cdots \widetilde{B_{\lambda,\mu_k}}^\fr$ --- achieving this upper bound. 
    
    By the induction hypothesis, for every $\nu\in\mathcal{N}$, $\widetilde{B_{\lambda,\nu}}^\fr$ is a polynomial of degree $\le k$ with integer coefficients in $\widetilde{B_{\lambda,1}}^\fr,\ldots,\widetilde{B_{\lambda,\nu_1}}^\fr$, and the latter functions are a subset of $\widetilde{B_{\lambda,1}}^\fr,\ldots,\widetilde{B_{\lambda,\mu_1}}^\fr$ as $\nu_1\le\mu_1$. Moreover, in these polynomials the sum $t_1+\ldots+t_j$ in every monomial is at most $|\nu|<|\mu|$. We conclude that the monomial $\widetilde{B_{\lambda,\mu_1}}^\fr \cdots \widetilde{B_{\lambda,\mu_k}}^\fr$ survives in the right hand side of \eqref{eq:expression for B_lambda,mu^fr}, and $\widetilde{B_{\lambda,\mu}}^\fr$ is given by a polynomial of degree exactly $k$ in $\widetilde{B_{\lambda,1}}^\fr, \ldots, \widetilde{B_{\lambda,\mu_1}}^\fr$, in which every monomial also satisfies the $(t_1+\ldots+t_j)$-condition.
\end{proof}

Given Lemma \ref{lem:poly degree k}, in order to establish Theorem \ref{thm:Nir}, it remains to prove that $(i)$ the polynomial $p_{\lambda,\mu}$ does not depend on $\lambda$ and that $(ii)$ the coefficients are polynomials in $q$. We first show the independence in $\lambda$.

\begin{lem} \label{lem:poly does not depend on lambda}
    Fixing $q$, the polynomial $p_{\lambda,\mu}$ does not depend on $\lambda$.
\end{lem}

\begin{proof}
We prove the lemma, again, by induction on $m=|\mu|$. The polynomial \eqref{eq:poly for B in R} we have for $\widetilde{B_{\lambda,\mu}}$ certainly does not depend on $\lambda$ (nor on $q$). By the induction hypothesis and Corollary \ref{cor:nu contained in mu}, for every $\nu\in\mathcal{N}=\mathcal{N}_{\lambda,\mu}$, the polynomial corresponding to $\widetilde{B_{\lambda,\nu}}^\fr$ does not depend on $\lambda$.  It thus remains to show that the multiset $\mathcal{N}=\mathcal{N}_{\lambda,\mu}$ of partitions does not depend on $\lambda$.

By the proof of Proposition \ref{prop:R and R^fr}, the multiset $\cal{N}$ is in one-to-one correspondence with the set of subspaces $0\ne\Omega\le \fq^m$ such that $\Omega B_{\lambda,\mu}=\Omega$. These are, equivalently, the subspaces $0\ne\Omega\le \fq^m$ such that $\Omega (B_{\lambda,\mu}-\lambda I)\le\Omega$. Thus, the set of valid subspaces $\Omega$ does not depend on $\lambda$. Because the matrices $S=S_\Omega$ and $T=T_\Omega$ from the proof of Proposition \ref{prop:R and R^fr} are determined by $\Omega$, and $ST=I_k$ (here $k=m-\dim(\Omega)$), we get that
\[
    C_\Omega-\lambda I_k = SB_{\lambda,\mu}T-\lambda I_k=S(B_{\lambda,\mu}-\lambda I_m)T.
\]
As the right hand side of the last equation does not depend on $\lambda$, we get that neither does the left hand side, so the partition $\nu$ corresponding to $C_\Omega$ does not depend on $\lambda$.
\end{proof}

As $p_{\lambda,\mu}$ does not depend on $\lambda$ (as long as $q$ is fixed), we now denote it, instead, by $^qp_\mu$. Similarly, we denote the multiset of partitions $\cal{N}_{\lambda,\mu}$ by $^q\cal{N}_\mu$.

The following lemma analyses the partition corresponding to the Jordan form of a nilpotent linear transformation on a finite dimensional vector space $V$. We assume that we are given the Jordan form in the action on an invariant subspace $W$ such that the action on $\nicefrac{V}{W}$ is given by a single Jordan block. We shall use it as a step in Lemma \ref{lem:poly has coefficients poly in q}.

\begin{lem} \label{lem:adding one generator to C}
    Assume that $T\colon V\to V$ is a nilpotent linear transformation on a finite dimensional vector space $V$, with an invariant subspace $W$ so that the quotient space $\nicefrac{V}{W}$ admits a $T$-cyclic vector: there exists a vector $v\in V$ such that $(v+W)_T=\nicefrac{V}{W}$ or, equivalently, $V=W+(v)_T$.

    Assume further that the Jordan-partition corresponding to the action of $T$ on $W$ is given by $\nu=(\nu_1,\ldots,\nu_r)$, that is, $W$ has a basis\footnote{To agree with the convention that our matrices $B$ act from the right on $V_m$, we also denote the action of $T$ on a vector $v$ by $vT$ -- this is simply the image of $v$ through $T$.} 
    \begin{equation} \label{eq:basis of W}
        w_1,w_1T,\ldots,w_1T^{\nu_1-1},w_2,w_2T,\ldots,w_2T^{\nu_2-1},\ldots,w_r,w_rT,\ldots,w_rT^{\nu_r-1}
    \end{equation}
    with $w_1T^{\nu_1}=\ldots=w_rT^{\nu_r}=0$.
    
    Then the Jordan-partition $\mu$ corresponding to the action of $T$ on $V$ depends only on $s\defi \dim\big(\nicefrac{V}{W}\big)$, and on the \textit{support} of the linear combination expressing $vT^s\in W$ in terms of the above basis elements of $W$.
\end{lem}

By the argument in the proof of Corollary \ref{cor:nu contained in mu}, $\nu$ is contained in $\mu$. The following proof can be made to show, moreover, that for every $j\in \mathbb{Z}_{\ge1}$ we have $\mu'_j\in\{\nu'_j,\nu'_j+1\}$. We do not use this additional property so we do not bother proving it explicitly.

\begin{proof}
    First, a basis for $V$ is given by the basis \eqref{eq:basis of W} of $W$ together with $v,vT,\ldots,vT^{s-1}$. By Lemma \ref{lem:extracting mu from action}, the columns of $\mu$ are the differences in the following sequence of dimensions
    \[
        \dim(V),\dim(VT),\dim(VT^2),\ldots
    \]
    Each term in this sequence is given by
    \begin{equation} \label{eq:VTj}
        \dim(VT^j)=\dim(WT^j)+\dim(\nicefrac{VT^j}{WT^j}).%=\nu'_j+\dim(\nicefrac{VT^j}{WT^j}).
    \end{equation}
    By Lemma \ref{lem:extracting mu from action}, the first summand is $\dim(WT^j)=|\nu|-\nu'_1-\ldots-\nu'_j$. To analyse the second summand, $\dim(\nicefrac{VT^j}{WT^j})$, we go over the basis elements $v,vT,\ldots,vT^{s-1}$ one by one, act by $T^j$ and see if each one is linearly dependent on the previous ones together with $WT^j$. Namely, we go over $vT^j,vT^{j+1},\ldots,vT^{j+s-1}$. 

    Now assume that for some $\ell\ge j$, $vT^\ell$ belongs to $WT^j+\text{span}\{vT^j,\ldots,vT^{\ell-1}\}$. We claim that $vT^\ell$ then belongs, in fact, to $WT^j$. Indeed, let
    \begin{equation} \label{eq:vT ell}
        vT^\ell = w + \alpha_1vT^{t_1}+\alpha_2vT^{t_2}+\ldots+\alpha_rvT^{t_r}
    \end{equation}    
    be a linear combination with smallest $r$ satisfying $w\in WT^j$ and $j\le t_1<t_2<\ldots<t_r\le \ell-1$. Assume towards contradiction that $r\ge1$. We may express each of the terms in either side of \eqref{eq:vT ell} in the basis we have for $V$. Notice first that $s\le t_1$ for otherwise $vT^{t_1}$ is a basis element of $V$ not appearing in any other term in \eqref{eq:vT ell}. So $vT^{t_1}\in W$ and is given by a linear combination in the basis elements \eqref{eq:basis of W} of $W$. Assume that $w_kT^p$ is an element in the support of this combination with $p$ minimal. By the minimality of $r$, $vT^{t_1}\notin WT^j$, which means that $p<j$. But then, $w_kT^p$ does not appear in the support of any other term in \eqref{eq:vT ell}: not in any $vT^t$ for $t\in \{t_2,t_3,\ldots,t_r,\ell\}$ where the minimal "power" of $T$ in the support is at least $p+t-t_1>p$, and not in $w\in WT^j$. This is a contradiction. 

    We conclude that $vT^\ell$ is linearly dependent in the former vectors and $WT^j$ if and only if it belongs to $WT^j$. Now whether or not 
    $vT^\ell \in WT^j$ is clearly determined only by the \textit{support} of the linear combination of the basis elements \eqref{eq:basis of W} of $W$ expressing $vT^s$. Hence for every $j$, $\dim(VT^j)$ is determined by that support, so $\mu$ is determined by it.
\end{proof}

The following lemma is the last step in proving Theorem \ref{thm:Nir}.

\begin{lem} \label{lem:poly has coefficients poly in q}
   The coefficients of $^qp_\mu$ are $\mathbb{Z}$-polynomials in $q$.
\end{lem}

\begin{proof}
    Denote $\mu=(\mu_1,\ldots,\mu_k)\vdash m$. We use, again, induction on $m$. The beginning of the argument is the same as in the proof of Lemma \ref{lem:poly does not depend on lambda}: in \eqref{eq:expression for B_lambda,mu^fr}, $\widetilde{B_{\lambda,\mu}}$ is given by a $\mathbb{Z}$-polynomial \eqref{eq:poly for B in R} which is independent of $q$, and by the induction hypothesis, for every $\nu\in\mathcal{N}_{\lambda,\mu}={^q}\cal{N}_\mu$, the polynomial $^qp_\nu$ corresponding to $\widetilde{B_{\lambda,\nu}}^\fr$ has coefficients which are $\mathbb{Z}$-polynomials in $q$.  It thus remains to show that in the multiset ${^q}\cal{N}_\mu$ of partitions, the multiplicity of every $\nu$ is given by a $\mathbb{Z}$-polynomial in $q$. We will use induction on $m$ also for this claim.

    Recall that any instance of some partition in ${^q}\cal{N}_\mu$ corresponds to some $\Omega\le V_m=\fq^m$ satisfying $\Omega=\Omega B_{\lambda,\mu}$, or, equivalently, $\Omega\ge\Omega B_\mu$, where 
    \[
        B_\mu\defi B_{\lambda,\mu}-\lambda I_m
    \]
    is the nilpotent $0,1$-matrix with zeros everywhere except for, possibly, in some entries in the subdiagonal. For example, if $\mu=(3,2)$, then
    \[
        B_\mu = \begin{pmatrix}
            0 & 0 & 0 & 0 & 0 \\
            1 & 0 & 0 & 0 & 0 \\
            0 & 1 & 0 & 0 & 0 \\
            0 & 0 & 0 & 0 & 0 \\
            0 & 0 & 0 & 1 & 0
        \end{pmatrix}.
    \]
    The subspace $\Omega$ corresponds to an instance of the partition $\nu\in {^q}\cal{N}_\mu$ if $\nu$ is the Jordan-partition corresponding to the action of $B_\mu$ on $\nicefrac{V_m}{\Omega}$.

    Next, let $\mu_{<k}\vdash m-\mu_k$ denote the partition obtained from $\mu$ by removing the last row, and let $W\le V_m$ be the $B_\mu$-invariant subspace corresponding to $B_{\mu_{<k}}$, namely, the one supported on the first $m-\mu_k$ coordinates in $V_m$. For every $\Omega\le V_m$ with $\Omega\ge \Omega B_\mu$, the subspace $\Omega\cap W$ is a $B_\mu$-invariant subspace of $W$. This subspace corresponds to an instance of some partition $\rho$ in ${^q}{\cal N}_{\mu_{<k}}$: the Jordan-partition corresponding to the action of $B_\mu$ on $W/(\Omega\cap W)$.
    
    Our strategy is to count the multiplicity of $\nu$ in ${^q}\cal{N}_\mu$, namely, the number of $\Omega\le V_m$ with $\Omega\ge \Omega B_\mu$ corresponding to $\nu$, by grouping them according to the partition $\rho$. We have
    \begin{eqnarray*}
        \text{multiplicity of }\nu\text{ in } {^q}{\cal N}_\mu &=& \# \{\Omega\le V_m\,\mid\, \Omega\ge\Omega B_\mu, \nu \text{ is the Jordan partition of } B_\mu \curvearrowright \nicefrac{V_m}{\Omega} \} \\
        &=& \sum_{\Delta\le W\colon\Delta\ge\Delta B_\mu} \#\left\{ \Omega\le V_m\,\biggr\vert\, 
        \begin{gathered}
            \Omega\ge\Omega B_\mu,\, \Omega\cap W=\Delta \\
            \nu \text{ is the Jordan partition of } B_\mu \curvearrowright \nicefrac{V_m}{\Omega}
        \end{gathered} \right\}.
    \end{eqnarray*}
    We will show that $(i)$ the number of $\Omega$ as above corresponding to some $\Delta$ is a $\mathbb{Z}$-polynomial in $q$, and that $(ii)$ this polynomial depends only on the partition $\rho$ corresponding to $\Delta$ and not on $\Delta$ itself. This would suffice as, by induction, the multiplicity of $\rho$ in ${^q}{\cal N}_{\mu_{<k}}$ is given by a $\mathbb{Z}$-polynomial in $q$, and the set of possible $\rho$'s is finite and independent of $q$ (by Corollary \ref{cor:nu contained in mu}, $\rho$ is contained in $\mu_{<k}$).

    Indeed, for any $\Omega\le V_m$ with $\Omega\ge \Omega B_\mu$, as $V_m=W+(e_m)_{B_\mu}$ (here $e_m$ is the last standard unit vector in $V_m$), we also have 
    \[
        \nicefrac{V_m}{\Omega} = \nicefrac{W}{\Omega\cap W}+{(\overline{e_m})_{B_\mu}},
    \]
    where $\overline{e_m}\defi e_m+\Omega$ is the coset of $e_m$ in $\nicefrac{V_m}{\Omega}$.

    But now we are exactly in the situation described in Lemma \ref{lem:adding one generator to C}: assuming we know the partition $\rho$ corresponding to the action of $B_\mu$ on $\nicefrac{W}{\Omega\cap W}$ and we fix a nice linear basis for $\nicefrac{W}{\Omega\cap W}$ as in that lemma, we have that $\nu$, the partition corresponding to the action of $B_\mu$ on $\nicefrac{V_m}{\Omega}$ is completely determined by the following data: $(i)$ the smallest $s$ with $\overline{e_m}B_\mu^{~s}\in \nicefrac{W}{\Omega\cap W}$ (this is, equivalently, $s=\dim(\nicefrac{V_m}{\Omega})-\dim(\nicefrac{W}{\Omega\cap W})$, and $(ii)$ the \textit{support} of the linear combination expressing $\overline{e_m}B_\mu^s$ in the above basis of $\nicefrac{W}{\Omega\cap W}$. Clearly, the possibilities for these data are finite and independent of $q$. Moreover,  for each possibility, the number of possible linear combinations with a given support is a polynomial in $q$ (in fact, it is precisely $(q-1)$ to the size of the support), and the linear combination completely determines $\Omega$. This completes the proof of the lemma, and also of Theorem \ref{thm:Nir}.
\end{proof}

\begin{remark}
    We stress that in the final paragraph of the proof of Lemma \ref{lem:poly has coefficients poly in q}, when we consider the linear combination of the basis elements of 
    $\nicefrac{W}{\Omega\cap W}$ expressing $\overline{e_m}B_\mu^s$, not every support is valid. As $e_mB_\mu^sB_\mu^{\mu_k-s}=e_mB_\mu^{\mu_k}=0$, every basis element $w\in\nicefrac{W}{\Omega\cap W}$ in the support must satisfy $wB_\mu^{\mu_k-s}=0$.
\end{remark}

%%%%%%%%%%%%%%%%%%%%%%%%%%%%%%%%%%%%%%%%%%%%%%%%%%%%%%%%%%%%%
% Appendix: An algebraic interpretation for \btil
%%%%%%%%%%%%%%%%%%%%%%%%%%%%%%%%%%%%%%%%%%%%%%%%%%%%%%%%%%%%%

\begin{appendices}

\section{The meaning of the generators of $\R$} \label{sec:btil is size of kernel}

In this short Appendix we record some illuminating  features of the basis elements $\btil$ of $\R$. They are mere observations, not relying on the main results of this paper.

To state these features, we want some canonical form for the elements of $G_m$. As $\mathbb{F}_q$ is not algebraically closed, the ordinary Jordan normal form is not always available. However, the \textbf{generalized} Jordan canonical form \cite{robinson1970generalized} is. For any monic $f=x^d+a_{d-1}x^{d-1}\ldots+a_0\in \mathbb{F}_q[x]$, the companion matrix $C_f$ of $f$ is
\[
    C_f=\begin{pmatrix}
        0      & 0      & \cdots & 0      & -a_0 \\
        1      & 0      & \cdots & 0      & -a_1 \\
        0      & 1      & \cdots & 0      & -a_2 \\
        \vdots & \vdots & \ddots & \ddots & \vdots \\
        0      & 0      & \cdots & 1      & -a_{d-1}
    \end{pmatrix}.
\]
If the polynomial $f$ is irreducible, a generalized Jordan block corresponding to $f$ is a block matrix of the form
\begin{equation} \label{eq:generalized Jordan block}
    \begin{pmatrix}
        C_f    & 0      & \cdots & 0 \\
        U      & C_f    & \cdots & 0 \\
        \vdots & \ddots & \ddots & \vdots \\
        0      & \cdots & U      & C_f
    \end{pmatrix},
\end{equation}
where $U$ is a matrix whose sole non-zero entry is a 1 in the upper right hand corner.
A matrix in generalized Jordan canonical form is a block diagonal matrix, where each block is a generalized Jordan block. By \cite{robinson1970generalized}, every $B\in G_m$ is conjugate to a matrix in generalized Jordan canonical form. 

\begin{prop} \label{prop:value of btil(g)}
    Let $B\in G_m$ and $g\in G_N$. Then $\btil_N(g)$ is an integral power of $q$. 
    
    If the generalized Jordan canonical form of $B$ is made of a single block\footnote{This is equivalent to that the characteristic polynomial of $B$ is equal to its minimal polynomial.} and $f\in \mathbb{F}_q[x]$ is the characteristic polynomial of $B$, then
    \begin{equation} \label{eq:btil as size of kernel}
        \btil_N(g)=\left|\ker(f(g))\right|.
    \end{equation}
\end{prop}

Of course, if the generalized Jordan form of $B$ consists of the generalized Jordan blocks $B_1,\ldots,B_r$, then by \eqref{eq:multiplication in R}, $\btil=\widetilde{B_1}\cdots\widetilde{B_r}$. In this sense, Proposition \ref{prop:value of btil(g)} leads to an interpretation of the value of $\btil_N(g)$ for arbitrary $B$.

\begin{proof}
    By definition, $\btil_N(g)$ is the size of the space of solutions to a system of linear equations $Mg=BM$ over $\mathbb{F}_q$ (with $mN$ variables), hence it is a power of $q$.

    Now, $\eqref{eq:btil as size of kernel}$ is equivalent to that the dimension of the space of solutions to $Mg=BM$ is equal to the dimension of $\ker(f(g))$ (both over $\mathbb{F}_q$). These two dimensions are unaltered if we extend the field, for the equations are defined over the base field $\mathbb{F}_q$ and $f$ has coefficients in $\mathbb{F}_q$. Hence it is enough to prove that these two dimensions are identical when we work over the algebraic closure $\overline{\mathbb{F}_q}$.

    Assume that $B$ has the form \eqref{eq:generalized Jordan block} with  $f\in\mathbb{F}_q[x]$ a monic irreducible polynomial of degree $d$, and $B\in G_{dr}$ (so $B$ is made of $r\times r$ blocks of size $d\times d$ each). Let $\lambda_1, \ldots, \lambda_d \in\overline{\mathbb{F}_q}$ be the roots of $f$ (note the roots are distinct as every irreducible polynomial over a finite field is separable.) Then, in $\gl_{dr}(\overline{\mathbb{F}_q})$, $B$ is conjugated to a matrix in ordinary Jordan normal form with $d$ distinct Jordan blocks, one block of size $r$ with eigenvalue $\lambda_i$ for every $i=1,\ldots,d$. Denote these blocks $B_1,\ldots,B_d$. By \eqref{eq:multiplication in R} and Lemma \ref{lem:conjugates give the same function}, $\btil_N(g)=\prod_{i=1}^d \widetilde{B_i}_N(g)$.    
    
    For every $g\in \gl_N(\overline{\mathbb{F}_q})$, $$\dim\ker(f(g))=\sum_{i=1}^d\dim\ker((g-\lambda_i)^r).$$ 
    So it is enough to prove the equality of dimensions for $$B=\begin{pmatrix}
        \lambda    & 0      & \cdots & 0 \\
        1      & \lambda    & \cdots & 0 \\
        \vdots & \ddots & \ddots & \vdots \\
        0      & \cdots & 1      & \lambda
    \end{pmatrix},$$ 
    a single Jordan block of size $r$ with some $\lambda\in \overline{\mathbb{F}_q}$.

    Let $v_1,\ldots,v_r$ be the rows of the matrix $M\in\textrm{Mat}_{r\times N}(\overline{\fq})$. The equation $Mg=BM$ translates to that 
    \[
        \begin{pmatrix}
            {-}~~v_1g~~{-}\\
            {-}~~v_2g~~{-} \\
            \vdots \\
            {-}~~v_rg~~{-}
        \end{pmatrix} = 
        \begin{pmatrix}
            {-}~~\lambda v_1~~{-}\\
            {-}~~v_1+\lambda v_2~~{-} \\
            \vdots \\
            {-}~~v_{r-1}+\lambda v_r~~{-}
        \end{pmatrix},
    \]
    so $v_{r-1}=v_r(g-\lambda), \ldots, v_1=v_2(g-\lambda)$ and $0=v_1(g-\lambda)$. So in every solution to the equations, $v_1,\ldots,v_{r-1}$ are determined by $v_r$, and $v_r$ defines a valid solution if and only if $v_r\in \ker\left((g-\lambda)^r\right)$.
\end{proof}

We remark that the interpretation above can also be given a more combinatorial flavor. For example, if $B$ is a single, ordinary Jordan block of size $r$ with some $\lambda\in\fq^*$, and $g$ has precisely $k$ Jordan block with eigenvalue $\lambda$ and sizes $a_1,\ldots,a_k$, then 
\[
    \btil_N(g)=|\ker((g-\lambda)^r)|=q^{\sum_{i=1}^k\min(r,a_i)}.
\]

\end{appendices}

\bibliographystyle{alpha}
\bibliography{stable_char_paper}

\noindent Danielle Ernst-West, School of Mathematical Sciences, Tel
Aviv University, Tel Aviv, 6997801, Israel\\
\texttt{daniellewest@mail.tau.ac.il }~\\

\noindent Doron Puder, School of Mathematical Sciences, Tel Aviv University,
Tel Aviv, 6997801, Israel\\
and School of Mathematics, Institute for Advanced Study, Princeton, NJ 08540, USA \\
\texttt{doronpuder@gmail.com}~\\

\noindent Yotam Shomroni, School of Mathematical Sciences, Tel Aviv
University, Tel Aviv, 6997801, Israel\\
\texttt{yotam.shomroni@gmail.com}
\end{document}